\documentclass[11pt]{article}

\usepackage{amsmath,diagbox}
\usepackage{indentfirst}
\usepackage{mathrsfs}
\usepackage{amsfonts}
\usepackage{arydshln}
\usepackage{enumerate}
\usepackage{setspace}
\usepackage{amssymb,amsthm,cases}
\usepackage{amsbsy,pifont}
\usepackage{latexsym,diagbox,tabularx}
\usepackage{amsmath,amscd,bbm,multirow}
\usepackage{graphicx,epsfig,extarrows,dsfont,mathtools}
\usepackage[final,allcolors=blue,colorlinks=true]{hyperref}
\usepackage[normalem]{ulem}
\input xy
\usepackage{caption,subcaption} 

\usepackage{enumerate}

\numberwithin{equation}{section}

\theoremstyle{plain}
\newtheorem{theorem}{Theorem}[section]

\newtheorem{corollary}[theorem]{Corollary}

\newtheorem{lemma}[theorem]{Lemma}
\newtheorem{proposition}[theorem]{Proposition}

\theoremstyle{definition}
\newtheorem{definition}[theorem]{Definition}

\newtheorem{remark}[theorem]{Remark}

\DeclareSymbolFont{txfontsA}{U}{txmia}{m}{it}
\SetSymbolFont{txfontsA}{bold}{U}{txmia}{bx}{it}
\DeclareFontSubstitution{U}{txmia}{m}{it}
\DeclareMathSymbol{\upalpha}{\mathord}{txfontsA}{"0B}
\DeclareMathSymbol{\upbeta}{\mathord}{txfontsA}{"0C}
\DeclareMathSymbol{\upgamma}{\mathord}{txfontsA}{"0D}
\DeclareMathSymbol{\updelta}{\mathord}{txfontsA}{"0E}
\DeclareMathSymbol{\upepsilon}{\mathord}{txfontsA}{"0F}
\DeclareMathSymbol{\upzeta}{\mathord}{txfontsA}{"10}
\DeclareMathSymbol{\upeta}{\mathord}{txfontsA}{"11}
\DeclareMathSymbol{\uptheta}{\mathord}{txfontsA}{"12}
\DeclareMathSymbol{\upiota}{\mathord}{txfontsA}{"13}
\DeclareMathSymbol{\upkappa}{\mathord}{txfontsA}{"14}
\DeclareMathSymbol{\uplambda}{\mathord}{txfontsA}{"15}
\DeclareMathSymbol{\upmu}{\mathord}{txfontsA}{"16}
\DeclareMathSymbol{\upnu}{\mathord}{txfontsA}{"17}
\DeclareMathSymbol{\upxi}{\mathord}{txfontsA}{"18}
\DeclareMathSymbol{\uppi}{\mathord}{txfontsA}{"19}
\DeclareMathSymbol{\uprho}{\mathord}{txfontsA}{"1A}
\DeclareMathSymbol{\upsigma}{\mathord}{txfontsA}{"1B}
\DeclareMathSymbol{\uptau}{\mathord}{txfontsA}{"1C}
\DeclareMathSymbol{\upupsilon}{\mathord}{txfontsA}{"1D}
\DeclareMathSymbol{\upphi}{\mathord}{txfontsA}{"1E}
\DeclareMathSymbol{\upchi}{\mathord}{txfontsA}{"1F}
\DeclareMathSymbol{\uppsi}{\mathord}{txfontsA}{"20}
\DeclareMathSymbol{\upomega}{\mathord}{txfontsA}{"21}
\DeclareMathSymbol{\upvarepsilon}{\mathord}{txfontsA}{"22}
\DeclareMathSymbol{\upvartheta}{\mathord}{txfontsA}{"23}
\DeclareMathSymbol{\upvarpi}{\mathord}{txfontsA}{"24}
\DeclareMathSymbol{\upvarrho}{\mathord}{txfontsA}{"25}
\DeclareMathSymbol{\upvarsigma}{\mathord}{txfontsA}{"26}
\DeclareMathSymbol{\upvarphi}{\mathord}{txfontsA}{"27}                   

\DeclareSymbolFont{ugmL}{OMX}{mdugm}{m}{n}
\SetSymbolFont{ugmL}{bold}{OMX}{mdugm}{b}{n}
\DeclareMathAccent{\wideparen}{\mathord}{ugmL}{"F3}

\def\pt{\partial}

\def\ra{\rightarrow}

\def\s{\subseteq}

\def\e{\epsilon}

\def\ol{\overline}

\def\vp{\varphi}

\def\pt{\partial}

\def\om{\omega}
\def\Om{\Omega}
\def\la{\lambda}
\def\al{\alpha}
\def\be{\beta}
\def\de{\delta}
\def\ga{\gamma}

\def\Ga{\Gamma}
\def\La{\Lambda}

\def\ts{\times}

\def\iy{\infty}

\def\f{\frac}

\def\df{\mathrm d}

\def\hra{\hookrightarrow}

\def\mcA{\mathcal{A}}

\def\mcR{\mathcal{R}}

\def\mcE{\mathcal{E}}
\def\mcD{\mathcal{D}}

	\DeclareMathOperator{\Div}{div}

	\DeclareMathOperator{\Ran}{Ran}

	\newcommand{\R}{\mathbb R}
	\newcommand{\N}{\mathbb N}

\begin{document}

\title{Null controllability of a class of fully degenerate parabolic equations on the cube}

\author{ Dong-Hui Yang$^a$, Mengze Gu$^b$\footnote{\small
		The corresponding author. Email: gu3059268737@gmail.com}, Yuanzhi Zhou$^a$\\
$^a${\it School of Mathematics and Statistics, Central South University}\\
{\it Changsha 410075, P.R. China}\\
$^b${\it School of Mathematics and Science, Changsha Normal University}\\ {\it  Changsha 410083, China}}

\date{}

\maketitle{}

\begin{abstract}
	In this paper, we investigate the internal null controllability of a class of two-dimensional degenerate parabolic equations on the unit square $\Omega = (0,1)^2$, governed by the operator $\partial_t \varphi - \operatorname{div}(x_2^\alpha \nabla \varphi)$ with $\alpha \in (0,1)$ and homogeneous Dirichlet boundary conditions. Unlike existing higher-dimensional studies where degeneracy occurs only in a single entry of the diffusion matrix, all diagonal components of the diffusion coefficient matrix here degenerate simultaneously along the boundary segment $(0,1) \times \{0\}$. To overcome the analytical obstacles caused by this full degeneracy, we combine the method of separation of variables with a new one-dimensional parameterized Carleman estimate that explicitly tracks the dependence on the transverse eigenvalues. By integrating this modal observability estimate with a one-dimensional spectral inequality via the Lebeau-Robbiano strategy, we establish the exact null controllability of the system in any positive time $T > 0$ with a localized control acting on an arbitrary non-empty open rectangle $\omega \s \Omega$.
	
	\medskip
	
	
	\medskip
	
	
\end{abstract}

\section{Introduction}

The controllability of partial differential equations has been an active branch of mathematics and engineering since the 1960s. The controllability of non-degenerate (or uniformly parabolic) partial differential equations has been extensively investigated in the literature; see, for instance, \cite{Barcena,Beauchard1,Coron,Dautray,Fursikov,Komornik,Lasiecka,Lasiecka1,Lions,Lu,Miller,Privat,Rousseau,Skrepek,Weiss,Yang5,Yao,Zuazua}. Since the 2000s, significant progress has been made on the controllability of degenerate partial differential equations. In particular, there is a vast body of work on the one-dimensional case \cite{Akil,Alabau,Bai,Buffe,Cannarsa1,Fragnelli,Gueye,Guo,Zhang}, alongside several contributions to the two-dimensional case \cite{Araruna,Beauchard,Cannarsa,Yang,Yang1,Yang2,Yang4,Yang9,Yang10} and higher-dimensional cases \cite{Yang3,Yang7,Yang8}. J.-L. Lions' Hilbert Uniqueness Method establishes that exact null controllability is equivalent to a suitable observability inequality for the associated adjoint system, a principle that remains valid in the degenerate setting \cite{Lions,Yang1}.

However, all existing results in higher dimensions ($N\geq 2$) focus on equations where the degeneracy occurs only in a single component (one-tuple) of the diffusion matrix, such as \cite{Araruna,Beauchard,Cannarsa,Yang,Yang2,Yang4,Yang6,Yang8,Yang9}. To the best of our knowledge, the controllability of degenerate equations where all components (all-tuples) of the diffusion matrix degenerate simultaneously—such as model \eqref{11.15.1} below-remains completely open. Addressing this gap is the primary motivation of the present paper.

In this paper, we investigate the following degenerate parabolic equation:
\begin{equation}\label{11.15.1}
	\begin{cases}
		\partial_t \varphi - \Div(w\nabla \varphi) = f & \text{in } Q, \\
		\varphi = 0 & \text{on } \Sigma, \\
		\varphi(0) = \varphi^0 & \text{in } \Omega,
	\end{cases}
\end{equation}
where $\Omega = (0,1)^2$, $Q = \Omega \times (0,T)$ for some $T > 0$, and $\Sigma = \partial \Omega \times (0,T)$. Here, the boundary components are defined by
\begin{equation}\label{11.15.2}
	\Gamma_1^i = \{i\} \times (0,1)  \text{ and }   \Gamma_2^j = (0,1) \times \{j\}   \text{ for } i,j = 0,1,
\end{equation}
with $\Lambda = (0,1)$ and $\Gamma = \partial\Omega - \Gamma_2^0$. The weight function $w$ and the associated elliptic operator $\mcA$ are given by
\begin{equation}\label{11.15.3}
	w(x) = x_2^\alpha   \text{ and }  \mcA \varphi = -\Div(w\nabla\varphi),
\end{equation}
where $\alpha \in (0,1)$ is a given degeneracy parameter, $\varphi^0 \in L^2(\Omega)$ is the initial datum, and $f \in L^2(Q)$ is a given source term. Notice that the diffusion coefficient $w(x) = x_2^\alpha$ degenerates along the boundary portion $\Gamma_2^0$.

The remainder of this paper is organized as follows. In Section \ref{S2}, we present some preliminary results, including weighted Sobolev spaces, the method of separation of variables, and the existence and uniqueness of mild solutions to \eqref{11.15.1}. In Section \ref{S3}, we establish the parameterized Carleman estimate and the null controllability result for \eqref{11.15.1}.

The main results of this paper are Lemmas \ref{11.15.L3} and \ref{11.15.L4}, Corollary \ref{11.27.C1}, and Theorem \ref{11.27.T2}.

\section{Preliminaries}\label{S2}

In this section, we present some preliminary results, including weighted Sobolev spaces, the method of separation of variables, and the existence and uniqueness of mild solutions to \eqref{11.15.1}.

\subsection{Weighted Sobolev spaces}

Define 
\begin{equation*}
	H^1(\Om;w)=\left\{u\in L^2(\Om) \cap H^1_{\mathrm{loc}}(\Om)\colon \int_\Om w|\nabla u|^2\df x<+\iy\right\},
\end{equation*}
equipped with the inner product and norm
\begin{equation*}
	(u,v)_{H^1(\Om;w)}=\int_\Om \left(uv+w\nabla u\cdot \nabla v\right)\df x,\quad \|u\|_{H^1(\Om;w)}=(u,u)_{H^1(\Om;w)}^\f{1}{2}. 
\end{equation*}
Set 
\begin{equation*}
	\begin{split}
		H_0^1(\Om;w)
		&=\text{the closure of $C_0^\iy(\Om)$ in }H^1(\Om;w), \\ 
		H^{-1}(\Om;w)
		&=\text{the dual space of $H_0^1(\Om;w)$ with the pivot space } L^2(\Om).
	\end{split}
\end{equation*} 
Furthermore, we define the domain of the operator $\mcA$ by
\begin{equation*}
	D(\mcA)=\left\{u\in H_0^1(\Om;w)\colon \mcA u\in L^2(\Om) \text{ in the sense of distributions}\right\},
\end{equation*}
endowed with the graph norm $\|u\|_{D(\mcA)} = \left(\|u\|_{L^2(\Om)}^2 + \|\mcA u\|_{L^2(\Om)}^2\right)^{\f{1}{2}}$. It is easily verified that $H^1(\Om;w)$, $H_0^1(\Om;w)$, $H^{-1}(\Om;w)$, and $D(\mcA)$ are Hilbert spaces (see, e.g., \cite{GC,Heinonen}).

\begin{lemma}\label{11.15.L1}
	For each $u\in H_0^1(\Om;w)$, we have
	\begin{equation}\label{11.15.4}
		\int_\Om x_2^{\al-2}u^2\df x\leq \f{4}{(1-\al)^2}\int_\Om x_2^\al |\pt_{x_2}u|^2\df x\leq \f{4}{(1-\al)^2}\int_\Om w|\nabla u|^2\df x.
	\end{equation}
\end{lemma}

\begin{proof}
	By the density of $C_0^\iy(\Om)$ in $H_0^1(\Om;w)$ and Fatou's lemma (passing to an almost everywhere convergent subsequence), we only need to show \eqref{11.15.4} for $u\in C_0^\iy(\Om)$.
	
	For each $\be\in (\al,1)$, noting that $\al - \be < 0$, we have
	\begin{equation*}
		\int_s^1 x_2^{\al-\be-1}\df x_2=\left.\f{1}{\al-\be}x_2^{\al-\be}\right|_{x_2=s}^{x_2=1} = \f{1}{\be-\al}\left(s^{\al-\be} - 1\right) \leq \f{1}{\be-\al}s^{\al-\be}
	\end{equation*}
	for any $s \in (0,1)$. Since $u(x_1,0)=0$, using the fundamental theorem of calculus, H\"older's inequality, and Fubini's theorem, we compute
	\begin{equation*}
		\begin{split}
			\int_\Om x_2^{\al-2}u^2\df x &= \int_0^1\int_0^1 x_2^{\al-2}|u(x_1,x_2)-u(x_1,0)|^2\df x_2\df x_1\\
			&= \int_0^1\int_0^1 x_2^{\al-2} \left(\int_0^{x_2}\pt_{x_2} u(x_1,s)\df s\right)^2\df x_2\df x_1\\
			&\leq \int_0^1\int_0^1 x_2^{\al-2}\left(\int_0^{x_2} s^{\be}|\pt_{x_2}u(x_1,s)|^2\df s\right)\left(\int_0^{x_2} s^{-\be}\df s\right)\df x_2\df x_1\\
			&= \f{1}{1-\be}\int_0^1\int_0^1\int_0^{x_2} x_2^{\al-\be-1} s^\be|\pt_{x_2}u(x_1,s)|^2\df s\df x_2\df x_1\\
			&= \f{1}{1-\be}\int_0^1\int_0^1\int_s^{1} x_2^{\al-\be-1}s^\be|\pt_{x_2}u(x_1,s)|^2\df x_2\df s\df x_1\\
			&\leq \f{1}{(1-\be)(\be-\al)}\int_0^1\int_0^1 s^\al |\pt_{x_2}u(x_1,s)|^2\df s\df x_1.
		\end{split}
	\end{equation*}
	Choosing the optimal parameter $\be = \f{1+\al}{2} \in (\al,1)$, we obtain $(1-\be)(\be-\al) = \f{(1-\al)^2}{4}$. Substituting this back into the above inequality yields \eqref{11.15.4}.
\end{proof}

\begin{remark}\label{11.15.R1}
	Since $x_2 \in (0,1)$ and $\al \in (0,1)$, we have $1 \leq x_2^{\al-2}$ in $\Om$. By Lemma \ref{11.15.L1}, for each $u \in H_0^1(\Om;w)$, it holds that
	\begin{equation}\label{11.15.10}
		\int_\Om u^2 \df x \leq \f{4}{(1-\al)^2}\int_\Om x_2^\al |\pt_{x_2}u|^2 \df x \leq \f{4}{(1-\al)^2}\int_\Om w|\nabla u|^2 \df x. 
	\end{equation}
	This is the weighted Poincar\'e inequality in $H_0^1(\Om;w)$. Consequently, the seminorm $\left(\int_\Om w|\nabla u|^2 \df x\right)^{\f{1}{2}}$ is equivalent to the full norm $\|\cdot\|_{H^1(\Om;w)}$ on $H_0^1(\Om;w)$. Hereafter, we equip $H_0^1(\Om;w)$ with the equivalent norm
	\begin{equation*}
		\|u\|_{H_0^1(\Om;w)} := \left(\int_\Om w|\nabla u|^2 \df x\right)^{\f{1}{2}}.
	\end{equation*}
\end{remark}

\begin{lemma}\label{11.15.L2}
	The embedding $H_0^1(\Om;w)\hra L^2(\Om)$ is compact. 
\end{lemma}

\begin{proof}
	First, since $x_2 \in (0,1)$ and $\al \in (0,1)$, we have $1 \leq x_2^{\al-2}$ in $\Om$. By Lemma \ref{11.15.L1}, it follows that
	\begin{equation*}
		\|u\|_{L^2(\Om)}^2 \leq \int_\Om x_2^{\al-2} u^2 \df x \leq \f{4}{(1-\al)^2} \|u\|_{H^1(\Om;w)}^2 \text{ for all } u \in H_0^1(\Om;w),
	\end{equation*}
	which implies that the embedding $H_0^1(\Om;w) \hra L^2(\Om)$ is continuous.
	
	To prove compactness, let $\{u_n\}_{n=1}^\iy$ be a bounded sequence in $H_0^1(\Om;w)$, so that there exists a constant $M > 0$ satisfying $\|u_n\|_{H^1(\Om;w)} \leq M$ for all $n \geq 1$. For any $\de \in (0,1)$, set $\Om_\de = (0,1) \times (\de, 1)$ and $\Om_\de^c = \Om \setminus \Om_\de = (0,1) \times (0,\de]$. In $\Om_\de$, the weight satisfies $w(x) = x_2^\al \geq \de^\al > 0$, and hence
	\begin{equation*}
		\|u_n\|_{H^1(\Om_\de)}^2 = \int_{\Om_\de} \left(u_n^2 + |\nabla u_n|^2\right) \df x \leq \max\{1, \de^{-\al}\} \|u_n\|_{H^1(\Om;w)}^2 \leq \max\{1, \de^{-\al}\} M^2.
	\end{equation*}
	Thus, $\{u_n\}_{n=1}^\iy$ is bounded in $H^1(\Om_\de)$. By the classical Rellich-Kondrachov theorem, the embedding $H^1(\Om_\de) \hra L^2(\Om_\de)$ is compact. Applying a standard diagonalization argument over a sequence $\de_k = \f{1}{k+1}$ ($k = 1,2,\dots$), we can extract a subsequence $\{u_{n_j}\}_{j=1}^\iy$ that is a Cauchy sequence in $L^2(\Om_\de)$ for every $\de \in (0,1)$.
	
	On the other hand, in $\Om_\de^c$, we have $1 \leq \de^{2-\al} x_2^{\al-2}$ since $2 - \al > 0$. Applying Lemma \ref{11.15.L1} again gives
	\begin{equation*}
		\int_{\Om_\de^c} |u_{n_j} - u_{n_l}|^2 \df x \leq \de^{2-\al} \int_\Om x_2^{\al-2} |u_{n_j} - u_{n_l}|^2 \df x \leq \f{4\de^{2-\al}}{(1-\al)^2} \|u_{n_j} - u_{n_l}\|_{H^1(\Om;w)}^2 \leq \f{16 M^2}{(1-\al)^2} \de^{2-\al}
	\end{equation*}
	for all $j, l \geq 1$.
	
	Now, given any $\e  > 0$, since $2 - \al > 0$, we can choose $\de \in (0,1)$ sufficiently small such that
	\begin{equation*}
		\f{16 M^2}{(1-\al)^2} \de^{2-\al} < \f{\e ^2}{2}.
	\end{equation*}
	With this $\de$ fixed, since $\{u_{n_j}\}_{j=1}^\iy$ is Cauchy in $L^2(\Om_\de)$, there exists an integer $N \geq 1$ such that for all $j, l \geq N$,
	\begin{equation*}
		\int_{\Om_\de} |u_{n_j} - u_{n_l}|^2 \df x < \f{\e ^2}{2}.
	\end{equation*}
	Combining the above two estimates yields
	\begin{equation*}
		\|u_{n_j} - u_{n_l}\|_{L^2(\Om)}^2 = \int_{\Om_\de} |u_{n_j} - u_{n_l}|^2 \df x + \int_{\Om_\de^c} |u_{n_j} - u_{n_l}|^2 \df x < \e ^2, \text{ for all } j, l \geq N.
	\end{equation*}
	Therefore, $\{u_{n_j}\}_{j=1}^\iy$ is a Cauchy sequence in $L^2(\Om)$ and hence converges strongly in $L^2(\Om)$. This completes the proof.
\end{proof}

\begin{lemma}\label{11.16.L1}
	The embedding $H_0^1(\Om;w) \hra W_0^{1,1}(\Om)$ is continuous.
\end{lemma}

\begin{proof}
	For any $u \in C_0^\iy(\Om)$ and each $i \in \{1, 2\}$, applying the Cauchy-Schwarz inequality yields
	\begin{equation*}
		\begin{split}
			\int_\Om |\pt_{x_i} u|\df x 
			&\le \left(\int_\Om x_2^{-\al}\df x\right)^{\f{1}{2}} \left(\int_\Om x_2^\al |\pt_{x_i} u|^2\df x\right)^{\f{1}{2}} \le \f{1}{\sqrt{1-\al}} \|u\|_{H_0^1(\Om;w)},
		\end{split}
	\end{equation*}
where we used $\int_\Omega x_2^{-\alpha}\,\mathrm dx=\int_0^1 x_2^{-\alpha}\,\mathrm dx_2=\f{1}{1-\al}$ since $\alpha\in(0,1)$.  Moreover, by the standard Poincar\'e inequality on $W_0^{1,1}(\Om)$, we also have $\|u\|_{L^1(\Om)} \le C \|\nabla u\|_{L^1(\Om)} \le C \|u\|_{H_0^1(\Om;w)}$. The desired continuous embedding then follows directly from the density of $C_0^\iy(\Om)$ in $H_0^1(\Om;w)$ and $W_0^{1,1}(\Om)$.
\end{proof}

\subsection{Method of separation of variables}

In this subsection, we apply the method of separation of variables to determine the eigenfunctions and eigenvalues of the degenerate elliptic operator $\mcA$ in $L^2(\Om)$, which will be used to represent the solution of \eqref{11.15.1}. Let $\La = (0,1)$. 

Consider the two-dimensional eigenvalue problem
\begin{equation}\label{11.15.5}
	\begin{cases}
		-\Div(w\nabla \Phi) = \la \Phi, & \text{in } \Om, \\
		\Phi = 0, & \text{on } \pt\Om.
	\end{cases}
\end{equation}
Looking for non-trivial solutions of the separated form $\Phi(x_1,x_2) = X(x_1)Y(x_2)$, we substitute this ansatz into \eqref{11.15.5} and divide by $x_2^\al X(x_1)Y(x_2)$ to obtain
\begin{equation*}
	-\f{X''(x_1)}{X(x_1)} = \f{(x_2^\al Y'(x_2))' + \la Y(x_2)}{x_2^\al Y(x_2)} = \mu,
\end{equation*}
where $\mu \in \mathbb{R}$ is a separation constant. This splits \eqref{11.15.5} into two one-dimensional eigenvalue problems.

First, in the $x_1$-direction, $X(x_1)$ satisfies the classical Dirichlet Sturm--Liouville problem
\begin{equation*}
	\begin{cases}
		-X''(x_1) = \mu X(x_1), & \text{in } \La, \\
		X(0) = X(1) = 0.
	\end{cases}
\end{equation*}
Its eigenvalues and corresponding $L^2(\La)$-normalized eigenfunctions are explicitly given by
\begin{equation}\label{11.15.6}
	\mu_m = (m\pi)^2, \quad X_m(x_1) = \sqrt{2}\sin(m\pi x_1), \quad m \in \N^*,
\end{equation}
and $\{X_m\}_{m=1}^\iy$ forms an orthonormal basis of $L^2(\La)$.

Second, for each fixed $m \in \N^*$, substituting $\mu = \mu_m = (m\pi)^2$ into the $x_2$-equation leads to the one-dimensional degenerate eigenvalue problem
\begin{equation}\label{11.15.7}
	\begin{cases}
		-\left(w Y'(x_2)\right)' + \mu_m w Y(x_2) = \la Y(x_2), & \text{in } \La, \\
		Y(0) = Y(1) = 0.
	\end{cases}
\end{equation}
Define the one-dimensional weighted Sobolev space
\begin{equation*}
	H^1(\La;w) = \left\{u\in L^2(\La) \cap H^1_{\mathrm{loc}}(\La)\colon \int_\La w|\pt_{x_2}u|^2\df x_2<+\iy\right\},
\end{equation*}
equipped with the inner product and norm
\begin{equation*}
	(u,v)_{H^1(\La;w)} = \int_\La \big(uv+w(\pt_{x_2}u)(\pt_{x_2}v)\big)\df x_2, \quad \|u\|_{H^1(\La;w)} = (u,u)_{H^1(\La;w)}^\f{1}{2}. 
\end{equation*}
Set 
\begin{equation*}
	\begin{split} 
		H_0^1(\La; w) 
		&= \text{the closure of } C_0^\iy(\La) \text{ in } H^1(\La;w),\\
		H^{-1}(\La;w)
		&= \text{the dual space of } H_0^1(\La;w) \text{ with the pivot space } L^2(\La).
	\end{split} 
\end{equation*}
By a one-dimensional version of Lemmas \ref{11.15.L1} and \ref{11.15.L2}, the embedding $H_0^1(\La; w) \hra L^2(\La)$ is compact. Moreover, in view of the one-dimensional counterpart of Remark \ref{11.15.R1}, we equip $H_0^1(\La;w)$ with the equivalent norm
\begin{equation*}
	\|u\|_{H_0^1(\La;w)} := \left(\int_\La w|\pt_{x_2}u|^2\df x_2\right)^\f{1}{2}.
\end{equation*}
Consequently, for each fixed $m \in \N^*$, the associated symmetric, coercive bilinear form
\begin{equation*}
	a_m(Y, \widetilde{Y}) = \int_\La w \left(Y'\widetilde{Y}' + \mu_m Y\widetilde{Y}\right) \df x_2, \text{ for all } Y, \widetilde{Y} \in H_0^1(\La; w)
\end{equation*}
induces a positive self-adjoint operator on $L^2(\La)$ with compact resolvent. By the spectral theorem for compact self-adjoint operators, \eqref{11.15.7} admits an increasing sequence of simple positive eigenvalues $\{\la_{mn}\}_{n=1}^\iy$ satisfying
\begin{equation*}
	0 < \la_{m1} < \la_{m2} < \cdots < \la_{mn} < \cdots, \quad \la_{mn} \to +\iy \text{ as } n \to +\iy,
\end{equation*}
and the corresponding $L^2(\La)$-normalized eigenfunctions $\{Y_{mn}\}_{n=1}^\iy \s H_0^1(\La; w)$ form an orthonormal basis of $L^2(\La)$.

Finally, setting 
\begin{equation*}
	\Phi_{mn}(x_1, x_2) = X_m(x_1) Y_{mn}(x_2),  \mbox{ for }m, n \in \N^*, 
\end{equation*}
we conclude from standard tensor-product Hilbert space theory that $\{\Phi_{mn}\}_{m,n=1}^\iy \s D(\mcA)$ forms an orthonormal basis of $L^2(\Om)$ consisting of eigenfunctions of $\mcA$ associated with the eigenvalues $\{\la_{mn}\}_{m,n=1}^\iy$, that is, $\mcA \Phi_{mn} = \la_{mn} \Phi_{mn}$.

\begin{lemma}\label{11.15.L3}
	For each $m, n \in \N^*$, define
	\begin{equation}\label{11.15.8}
		\Phi_{mn}(x_1, x_2) = X_m(x_1) Y_{mn}(x_2), \quad (x_1, x_2) \in \Om.
	\end{equation}
	Then $\{\Phi_{mn}\}_{m,n=1}^\iy \s D(\mcA)$ forms an orthonormal basis of $L^2(\Om)$ as well as an orthogonal basis of $H_0^1(\Om;w)$, and each $\Phi_{mn}$ is an eigenfunction of $\mcA$ associated with the eigenvalue $\la_{mn}$, namely,
	\begin{equation}\label{11.15.9}
		\mcA \Phi_{mn} = \la_{mn} \Phi_{mn}  \text{ in } \Om.
	\end{equation}
\end{lemma}

\begin{proof}
	We divide the proof into four steps.
	
	{\it Step 1.} We prove $\Phi_{mn} \in D(\mcA)$ for all $m,n\in\N^*$.
	
	Fix $m, n \in \N^*$. Since $X_m \in H_0^1(\La)$ and $Y_{mn} \in H_0^1(\La; w)$, we obtain $\Phi_{mn} \in H_0^1(\Om; w)$. Next, for any test function $\phi \in C_0^\iy(\Om)$, using Fubini's theorem together with \eqref{11.15.6} and \eqref{11.15.7}, we obtain
	\begin{equation*}
		\langle \mcA \Phi_{mn}, \phi \rangle_{\mcD'(\Om), C_0^\iy(\Om)} 
		= \int_\Om w \nabla \Phi_{mn} \cdot \nabla \phi \df x = \la_{mn} \int_\Om \Phi_{mn}(x) \phi(x) \df x.
	\end{equation*}
	This shows that $\mcA \Phi_{mn} = \la_{mn} \Phi_{mn} \in L^2(\Om)$ in the sense of distributions. Therefore, $\Phi_{mn} \in D(\mcA)$ and \eqref{11.15.9} holds.
	
	{\it Step 2.} Orthonormality in $L^2(\Om)$.
	
	For any $m, n, k, l \in \N^*$, by Fubini's theorem and the orthonormality of $\{X_m\}_{m=1}^\iy$ and $\{Y_{mn}\}_{n=1}^\iy$ in $L^2(\La)$, we have
	\begin{equation*}
		(\Phi_{mn}, \Phi_{kl})_{L^2(\Om)} 
		= \left(\int_\La X_m(x_1) X_k(x_1) \df x_1\right) \left(\int_\La Y_{mn}(x_2) Y_{kl}(x_2) \df x_2\right) 
		= \de_{mk} \de_{nl},
	\end{equation*}
	where $\de_{ij}$ denotes the Kronecker delta, i.e., $\de_{mk}=1$ for $m=k$ and $\de_{mk}=0$ for $m\neq k$. Hence, $\{\Phi_{mn}\}_{m,n=1}^\iy$ is an orthonormal system in $L^2(\Om)$.
	
	{\it Step 3.} Completeness in $L^2(\Om)$.
	
	Suppose that $g \in L^2(\Om)$ satisfies $(g, \Phi_{mn})_{L^2(\Om)} = 0$ for all $m, n \in \N^*$. By Fubini's theorem, for each fixed $m \in \N^*$, the function
	\begin{equation*}
		G_m(x_2) := \int_\La g(x_1, x_2) X_m(x_1) \df x_1
	\end{equation*}
	belongs to $L^2(\La)$ and satisfies
	\begin{equation*}
		(G_m, Y_{mn})_{L^2(\La)} = \int_\La \left(\int_\La g(x_1, x_2) X_m(x_1) \df x_1\right) Y_{mn}(x_2) \df x_2 = (g, \Phi_{mn})_{L^2(\Om)} = 0
	\end{equation*}
	for all $n \in \N^*$. Since $\{Y_{mn}\}_{n=1}^\iy$ is an orthonormal basis of $L^2(\La)$, it follows that $G_m(x_2) = 0$ for a.e. $x_2 \in \La$ and every $m \in \N^*$. Since $\N^*$ is countable, there exists a null set $E \s \La$ such that for every $x_2 \in \La - E$, we have $g(\cdot, x_2) \in L^2(\La)$ and
	\begin{equation*}
		(g(\cdot, x_2), X_m)_{L^2(\La)} = G_m(x_2) = 0, \text{ for all } m \in \N^*.
	\end{equation*}
	Using the completeness of $\{X_m\}_{m=1}^\iy$ in $L^2(\La)$, we obtain $g(x_1, x_2) = 0$ for a.e. $x_1 \in \La$ and every $x_2 \in \La - E$. Consequently, $g = 0$ a.e. in $\Om$. Thus, $\{\Phi_{mn}\}_{m,n=1}^\iy$ is an orthonormal basis of $L^2(\Om)$.
	
	{\it Step 4.} Orthogonality and completeness in $H_0^1(\Om;w)$.
	
	Since $\Phi_{mn} \in D(\mcA)$ satisfies \eqref{11.15.9}, for any $m, n, k, l \in \N^*$ we have
	\begin{equation*}
		(\Phi_{mn}, \Phi_{kl})_{H_0^1(\Om;w)} 
		= (\mcA \Phi_{mn}, \Phi_{kl})_{L^2(\Om)} = \la_{mn} (\Phi_{mn}, \Phi_{kl})_{L^2(\Om)} 
		= \la_{mn} \de_{mk} \de_{nl},
	\end{equation*}
	which shows that $\{\Phi_{mn}\}_{m,n=1}^\iy$ is an orthogonal system in $H_0^1(\Om;w)$. Furthermore, if $u \in H_0^1(\Om;w)$ satisfies $(u, \Phi_{mn})_{H_0^1(\Om;w)} = 0$ for all $m, n \in \N^*$, then
	\begin{equation*}
		0 = (u, \Phi_{mn})_{H_0^1(\Om;w)} = \la_{mn} (u, \Phi_{mn})_{L^2(\Om)}, \text{ for all } m, n \in \N^*.
	\end{equation*}
	Since $\la_{mn} \geq \f{(1-\al)^2}{4} > 0$ by \eqref{11.15.10}, we obtain $(u, \Phi_{mn})_{L^2(\Om)} = 0$ for all $m, n \in \N^*$, which together with Step 3 implies $u = 0$. Therefore, $\{\Phi_{mn}\}_{m,n=1}^\iy$ is an orthogonal basis of $H_0^1(\Om;w)$. This completes the proof of the lemma.
\end{proof}

\begin{lemma}\label{11.15.L4}
	Let $\{\la_{mn}\}_{m,n=1}^\iy$ be the eigenvalues of \eqref{11.15.7}. Then there exist two positive constants $C_1, C_2 > 0$, depending only on $\al$, such that
	\begin{equation}\label{11.15.11}
		C_1 \left(m^{2-\al} + n^2\right) \leq \la_{mn} \leq C_2 \left(m^2 + n^2\right) \quad \text{for all } m, n \in \N^*.
	\end{equation}
\end{lemma}

\begin{proof}
	We divide the proof into three steps.
	
	{\it Step 1.} Explicit bounds for the unperturbed degenerate eigenvalues.
	
	Consider the auxiliary one-dimensional degenerate eigenvalue problem without the potential term $\mu_m w$:
	\begin{equation}\label{11.15.12}
		\begin{cases}
			-\left(x_2^\al \psi'(x_2)\right)' = \la_{0n} \psi(x_2), & \text{in } \La, \\
			\psi(0) = \psi(1) = 0,
		\end{cases}
	\end{equation}
	where the eigenvalues are arranged in increasing order $0 < \la_{01} < \la_{02} < \cdots < \la_{0n} \to +\iy$. Set
	\begin{equation*}
		\nu_\al := \f{1-\al}{2-\al} \in \left(0, \f{1}{2}\right), \quad \kappa_\al := \f{2-\al}{2} \in \left(\f{1}{2}, 1\right).
	\end{equation*}
	As established in \cite{Cannarsa1,Gueye}, the eigenvalues of \eqref{11.15.12} are explicitly given by
	\begin{equation*}
		\la_{0n} = \kappa_\al^2 j_{\nu_\al, n}^2 = \f{(2-\al)^2}{4} j_{\nu_\al, n}^2,  \text{ for all } n \in \N^*,
	\end{equation*}
	where $0 < j_{\nu_\al, 1} < j_{\nu_\al, 2} < \cdots < j_{\nu_\al, n} < \cdots$ denote the positive zeros of the Bessel function $J_{\nu_\al}$ of the first kind of order $\nu_\al$. Since $\nu_\al \in \left(0, \f{1}{2}\right)$, classical bounds on the zeros of Bessel functions (see \cite[Section 15.3]{Watson} or \cite{Cannarsa}) yield
	\begin{equation*}
		\left(n - \f{1}{4}\right)\pi < \left(n + \f{\nu_\al}{2} - \f{1}{4}\right)\pi \leq j_{\nu_\al, n} \leq n\pi,  \text{ for all } n \in \N^*.
	\end{equation*}
	Noting that $n - \f{1}{4} \geq \f{3}{4}n$ for every $n \in \N^*$, we deduce that
	\begin{equation}\label{11.15.13}
		c_1 n^2 \leq \la_{0,n} \leq c_2 n^2,  \text{ for all } n \in \N^*,
	\end{equation}
	where
	\begin{equation*}
		c_1 := \f{9(2-\al)^2 \pi^2}{64} > 0, \quad c_2 := \f{(2-\al)^2 \pi^2}{4} > 0.
	\end{equation*}
	
	{\it Step 2.} Comparison between $\la_{mn}$ and $\la_{0n}$ via the min--max principle.
	
	By the Courant--Fischer min--max characterization of eigenvalues for compact self-adjoint resolvents, we have
	\begin{equation*}
		\la_{0n} = \min_{\substack{V_n \s H_0^1(\La;w) \\ \dim V_n = n}} \max_{Y \in V_n - \{0\}} \f{\int_\La x_2^\al |Y'(x_2)|^2 \df x_2}{\int_\La |Y(x_2)|^2 \df x_2},
	\end{equation*}
	and similarly,
	\begin{equation*}
		\la_{mn} = \min_{\substack{V_n \s H_0^1(\La;w) \\ \dim V_n = n}} \max_{Y \in V_n - \{0\}} \f{\int_\La x_2^\al \left(|Y'(x_2)|^2 + \mu_m |Y(x_2)|^2\right) \df x_2}{\int_\La |Y(x_2)|^2 \df x_2},
	\end{equation*}
	where the minimum is taken over all $n$-dimensional subspaces $V_n$ of $H_0^1(\La;w)$. Since $0 < x_2^\al < 1$ for all $x_2 \in \La = (0,1)$ and $\mu_m = (m\pi)^2 > 0$, for every $Y \in H_0^1(\La;w) - \{0\}$ it holds that
	\begin{equation*}
		\f{\int_\La x_2^\al |Y'|^2 \df x_2}{\int_\La Y^2 \df x_2} 
		\leq \f{\int_\La x_2^\al \left(|Y'|^2 + \mu_m Y^2\right) \df x_2}{\int_\La Y^2 \df x_2} 
		\leq \f{\int_\La x_2^\al |Y'|^2 \df x_2}{\int_\La Y^2 \df x_2} + \mu_m.
	\end{equation*}
	Taking the maximum over $Y \in V_n - \{0\}$ and then the minimum over all $n$-dimensional subspaces $V_n \s H_0^1(\La;w)$, we obtain
	\begin{equation}\label{11.15.14}
		\la_{0n} \leq \la_{mn} \leq \la_{0n} + \mu_m, \text{ for all } m, n \in \N^*.
	\end{equation}
	Combining \eqref{11.15.13} and \eqref{11.15.14}, we immediately derive the upper bound
	\begin{equation}\label{11.15.15}
		\la_{mn} \leq c_2 n^2 + \pi^2 m^2 \leq C_2 \left(m^2 + n^2\right), \text{ with } C_2 := \max\{c_2, \pi^2\} = \pi^2,
	\end{equation}
	as well as the lower bound in terms of $n$:
	\begin{equation}\label{11.15.16}
		\la_{mn} \geq \la_{0n} \geq c_1 n^2, \text{ for all } m, n \in \N^*.
	\end{equation}
	
	{\it Step 3.} Lower bound in terms of $m$.
	
	Since the eigenvalues $\{\la_{mn}\}_{n=1}^\iy$ are strictly increasing with respect to $n$ for each fixed $m \in \N^*$, we have
	\begin{equation*}
		\la_{mn} \geq \la_{m1} = \min_{Y \in H_0^1(\La;w) - \{0\}} \f{\int_\La x_2^\al \left(|Y'|^2 + \mu_m Y^2\right) \df x_2}{\int_\La Y^2 \df x_2}.
	\end{equation*}
	Fix any $Y \in H_0^1(\La;w) - \{0\}$ and any parameter $\de \in (0,1)$, and decompose $\La = (0,\de) \cup [\de,1)$. 
	\begin{itemize}
		\item On the subinterval $[\de, 1)$, since $x_2 \geq \de$ and $\al > 0$, we have $x_2^\al \geq \de^\al$, and thus
		\begin{equation*}
			\int_\de^1 Y^2 \df x_2 \leq \f{1}{\mu_m \de^\al} \int_\de^1 x_2^\al \mu_m Y^2 \df x_2 \leq \f{1}{\mu_m \de^\al} \int_\La x_2^\al \mu_m Y^2 \df x_2.
		\end{equation*}
		\item On the subinterval $(0, \de)$, since $0 < x_2 < \de$ and $\al - 2 < 0$, we have $1 \leq \de^{2-\al} x_2^{\al-2}$. Applying the one-dimensional version of Lemma \ref{11.15.L1} yields
		\begin{equation*}
			\int_0^\de Y^2 \df x_2 \leq \de^{2-\al} \int_0^\de x_2^{\al-2} Y^2 \df x_2 \leq \de^{2-\al} \int_\La x_2^{\al-2} Y^2 \df x_2 \leq \f{4\de^{2-\al}}{(1-\al)^2} \int_\La x_2^\al |Y'|^2 \df x_2.
		\end{equation*}
	\end{itemize}
	Summing the two inequalities above, we find that
	\begin{align*}
		\int_\La Y^2 \df x_2 
		&= \int_0^\de Y^2 \df x_2 + \int_\de^1 Y^2 \df x_2\leq \f{4\de^{2-\al}}{(1-\al)^2} \int_\La x_2^\al |Y'|^2 \df x_2 + \f{1}{\mu_m \de^\al} \int_\La x_2^\al \mu_m Y^2 \df x_2 \\
		&\leq \max\left\{\f{4\de^{2-\al}}{(1-\al)^2}, \, \f{1}{\mu_m \de^\al}\right\} \int_\La x_2^\al \left(|Y'|^2 + \mu_m Y^2\right) \df x_2.
	\end{align*}
	Now choose $\de = (\mu_m)^{-\f{1}{2}} = \f{1}{m\pi}$. Since $m \geq 1$ and $\pi > 3$, we clearly have $\de \in (0,1)$, and the two exponents balance:
	\begin{equation*}
		\de^{2-\al} = (m\pi)^{-(2-\al)}, \quad \f{1}{\mu_m \de^\al} = (m\pi)^{-2} (m\pi)^\al = (m\pi)^{-(2-\al)}.
	\end{equation*}
	Consequently,
	\begin{equation*}
		\max\left\{\f{4\de^{2-\al}}{(1-\al)^2}, \, \f{1}{\mu_m \de^\al}\right\} = \max\left\{\f{4}{(1-\al)^2}, \, 1\right\} (m\pi)^{-(2-\al)} = \f{4}{(1-\al)^2} (m\pi)^{-(2-\al)},
	\end{equation*}
	where we used $\f{4}{(1-\al)^2} > 4 > 1$ for $\al \in (0,1)$. Substituting this back into the Rayleigh quotient gives
	\begin{equation}\label{11.15.17}
		\la_{mn} \geq \la_{m1} \geq \f{(1-\al)^2 \pi^{2-\al}}{4} m^{2-\al} =: c_3 m^{2-\al}, \text{ for all } m, n \in \N^*.
	\end{equation}
	
	Finally, adding \eqref{11.15.16} and \eqref{11.15.17} and dividing by $2$, we obtain
	\begin{equation*}
		\la_{mn} \geq \f{1}{2}\left(c_3 m^{2-\al} + c_1 n^2\right) \geq C_1 \left(m^{2-\al} + n^2\right), \text{ with } C_1 := \f{1}{2}\min\{c_1, c_3\} > 0.
	\end{equation*}
	Together with \eqref{11.15.15}, this establishes \eqref{11.15.11} and completes the proof.
\end{proof}

Let
\begin{equation}\label{11.26.4}
	\mathcal A^{(1)}u=-\partial_{x_1x_1}u,\quad\mathcal A^{(2)}u=-\partial_{x_2}(x_2^\alpha\partial_{x_2}u).
\end{equation}
The two-dimensional differential expression is\[\mathcal Au=-x_2^\alpha\partial_{x_1x_1}u-\partial_{x_2}(x_2^\alpha\partial_{x_2}u)=x_2^\alpha\mathcal A^{(1)}u+\mathcal A^{(2)}u,\]where the two expressions on the right act in their respective variables. As one-dimensional operators on $L^2(\Lambda)$, their domains are\[\begin{split}D(\mathcal A^{(1)})&=\{u\in H_0^1(\Lambda):\mathcal A^{(1)}u\in L^2(\Lambda)\},\\D(\mathcal A^{(2)})&=\{u\in H_0^1(\Lambda;w):\mathcal A^{(2)}u\in L^2(\Lambda)\}.\end{split}\]

Define the differential operators
\begin{equation}\label{11.27.1}
	\mcA_0^{(2)}u=-\pt_{x_2}(x_2^\al \pt_{x_2}u)  \text{ and }   \mcA_m^{(2)}u=-\pt_{x_2}(x_2^\al \pt_{x_2}u)+\mu_m x_2^\al u  \text{ for } m\in\N^*, 
\end{equation}
with their domains given by
\begin{equation*}
	D(\mcA_m^{(2)})=\left\{ u\in H_0^1(\La;w)\colon \mcA_m^{(2)}u\in L^2(\La)\right\}  \text{ for } m\in\N. 
\end{equation*}
Clearly, $D(\mcA_m^{(2)})=D(\mcA_0^{(2)})$ for all $m\in\N^*$.

\subsection{Mild solution}\label{S2.3}

Since $\mcA$ is a densely defined, self-adjoint, positive operator on $L^2(\Om)$ with domain $D(\mcA)$, it follows from standard semigroup theory (see, e.g., \cite{Pazy}) that $-\mcA$ generates a strongly continuous analytic contraction semigroup $\{S(t)\}_{t \geq 0}$ on $L^2(\Om)$. Moreover, in view of the orthonormal basis $\{\Phi_{mn}\}_{m,n=1}^\iy$ established in Lemma \ref{11.15.L3}, the semigroup $S(t)$ admits the explicit spectral representation
\begin{equation}\label{11.15.18}
	S(t)\psi = \sum_{m,n=1}^\iy e^{-\la_{mn} t} (\psi, \Phi_{mn})_{L^2(\Om)} \Phi_{mn}, \text{ for all } \psi \in L^2(\Om) \text{ and } t \geq 0.
\end{equation}

We now recall the definition of mild solutions to the inhomogeneous Cauchy problem \eqref{11.15.1}.

\begin{definition}\label{11.15.D_mild}
	Given $\vp^0 \in L^2(\Om)$ and $f \in L^1(0,T; L^2(\Om))$, a function $\vp \in C([0,T]; L^2(\Om))$ is called a \emph{mild solution} to \eqref{11.15.1} if it satisfies the variation-of-constants (Duhamel) formula
	\begin{equation}\label{11.15.19}
		\vp(t) = S(t)\vp^0 + \int_0^t S(t-s)f(s) \df s, \text{ for all } t \in [0,T].
	\end{equation}
\end{definition}

The well-posedness and basic energy regularity of the mild solution are summarized in the following proposition.

\begin{proposition}\label{11.15.P_mild}
	For any $\vp^0 \in L^2(\Om)$ and $f \in L^2(Q)$, the problem \eqref{11.15.1} admits a unique mild solution $\vp$ given by \eqref{11.15.19}, which also belongs to $L^2(0,T; H_0^1(\Om;w))$ and coincides with the unique weak solution of \eqref{11.15.1}. Furthermore, $\vp$ satisfies the energy estimate
	\begin{equation}\label{11.15.20}
		\max_{t \in [0,T]} \|\vp(t)\|_{L^2(\Om)}^2 + \|\vp\|_{L^2(0,T; H_0^1(\Om;w))}^2 \leq C_T \left( \|\vp^0\|_{L^2(\Om)}^2 + \|f\|_{L^2(Q)}^2 \right),
	\end{equation}
	where $C_T > 0$ is a constant depending only on $T$ (in fact, one may take $C_T = 2e^T$ or simply $C_T = 2 + T$ when dropping the lower-order Gronwall factor).
\end{proposition}

From Lemma \ref{11.15.L3}, the weak solution given by \eqref{11.15.19} admits the Fourier expansion
\begin{equation*}
	\vp(x_1,x_2,t)=\sum_{m\in\N^*} \vp_m(x_2,t)X_m(x_1),
\end{equation*}
where for each $m\in\N^*$, the Fourier coefficient $\vp_m(x_2,t)$ satisfies the one-dimensional degenerate parabolic problem
\begin{equation}\label{11.15.21}
	\begin{cases}
		\pt_t\vp_m-\pt_{x_2} \left(x_2^\al \pt_{x_2}\vp_m\right)+\mu_m x_2^\al \vp_m=f_m & \text{in }\La_T:=\La\ts (0,T),\\
		\vp_m(0,t)=\vp_m(1,t)=0 & \text{for } t\in (0,T),\\
		\vp_m(x_2,0)=\vp_m^0(x_2) & \text{in }\La, 
	\end{cases}
\end{equation}
with
\begin{equation*}
	f_m(x_2,t)=\left(f(\cdot,x_2,t),X_m\right)_{L^2(\La)},\quad \vp_m^0(x_2)=\left(\vp^0(\cdot,x_2),X_m\right)_{L^2(\La)}. 
\end{equation*}

\section{Carleman estimate and controllability}\label{S3}

In this section, we first establish a parameterized Carleman estimate for the modal problem \eqref{11.15.21}, and then deduce the null controllability of \eqref{11.15.1}.

Set
\begin{equation}\label{11.26.1}
	\be=\f{3-\al}{2},\quad q=\f{3\al-1}{2},\quad h=\f{1+\al}{2}. 
\end{equation}
Since $\al\in (0,1)$, we have $1<\be<2-\al$ and $q<\al$. Given $(a,b)\Subset (0,1)$, we choose $c, d$ such that $0<a<c<d<b<1$, and define
\begin{equation}\label{11.26.2}
	\eta(x_2)=
	\zeta (x_2) x_2^\be+(1-\zeta(x_2))(1-x_2)x_2^{-\al},
\end{equation}
where $\zeta=\zeta(x_2)\in C^\iy(\R)$ satisfies $0\leq \zeta\leq 1$,
\begin{equation*}
	\zeta=1 \text{ for } x_2\leq c, \text{ and }  \zeta=0 \text{ for } x_2\geq d. 
\end{equation*}
Then $\eta\in C^1([0,1])\cap C^\iy((0,1])$, $\eta>0$ in $(0,1)$, and $\eta(0)=\eta(1)=0$. 
Define 
\begin{equation}\label{11.26.3}
	\Theta(t)=\f{1}{[t(T-t)]^4},\quad \ga=\|\eta\|_{L^\iy(\La)}+1, \quad \xi(x_2,t)=\Theta(t)(\ga-\eta(x_2)). 
\end{equation}

For any $s > 0$, set
\begin{equation*}
	z=\vp_me^{-s\xi} \iff \vp_m=ze^{s\xi}. 
\end{equation*}
Then the first equation in \eqref{11.15.21} can be rewritten as
\begin{equation*}
	e^{-s\xi}f_m=P_1z+P_2z+Kz, 
\end{equation*}
where 
\begin{equation*}
	\begin{split}
		P_1z
		&=z_t-2s x_2^\al (\pt_{x_2}z)(\pt_{x_2}\xi)-sz\pt_{x_2}(x_2^\al \pt_{x_2}\xi),\\
		P_2z
		&=-\pt_{x_2}(x_2^\al \pt_{x_2}z)+sz\xi_t-s^2zx_2^\al(\pt_{x_2}\xi)^2, \\
		Kz
		&=\mu_mx_2^\al z. 
	\end{split}
\end{equation*}

\subsection{Justification of the integration by parts}\label{S3.1}

Before deriving the Carleman estimate, we first justify the validity of the integration by parts used in computing the inner products $(P_1 z, P_2 z)_{L^2(\La_T)}$ and $(P_1 z, Kz)_{L^2(\La_T)}$.

Let $\vp_m^0\in D(\mcA_m^{(2)})$ and $f_m\in H^1(0,T; H_0^1(\La;w))$. By standard regularity theory for linear parabolic equations (since $\mcA_m^{(2)}$ is a positive self-adjoint operator on $L^2(\La)$), the unique solution $\vp_m$ of \eqref{11.15.21} satisfies
\begin{equation*}
	\vp_m \in C([0,T]; D(\mcA_m^{(2)})) \cap H^1(0,T; H_0^1(\La;w)).
\end{equation*}
Since $e^{-s\xi}$ decays exponentially together with all its derivatives as $t \to 0^+$ and $t \to T^-$, the function $z = \vp_m e^{-s\xi}$ inherits the regularity
\begin{equation*}
	z \in C([0,T]; D(\mcA_m^{(2)})) \cap H^1(0,T; H_0^1(\La;w)),
\end{equation*}
with 
\begin{equation}\label{11.26.5}
	z(\cdot, 0) = z(\cdot, T) = 0  \text{ and }  z_t(\cdot, 0) = z_t(\cdot, T) = 0.
\end{equation}
Moreover, 
\begin{equation}\label{11.26.8}
	z=\pt_t z=0  \text{ on } \{0,1\}\ts (0,T). 
\end{equation}

Furthermore, using the definitions \eqref{11.26.1}-\eqref{11.26.2} of $\eta$, we have $|\pt_{x_2}\eta(x_2)| \leq C x_2^{\be-1} = C x_2^{\f{1-\al}{2}}$ and $|\pt_{x_2}(x_2^\al \pt_{x_2}\eta(x_2))| \leq C x_2^{\f{\al-1}{2}}$ for $x_2 \in (0,1)$. Combined with the Hardy-type inequality
\begin{equation}\label{11.26.6}
	\begin{split} 
		|z(x_2,t)|^2
		&=\left(\int_0^{x_2}\pt_{x_2}z(y,t)\df y\right)^2\\
		&\leq \left(\int_0^{x_2} y^\al |\pt_{x_2}z(y,t)|^2\df y\right) \left(\int_0^{x_2} y^{-\al}\df y\right)\leq \f{1}{1-\al}x_2^{1-\al}\|z(\cdot,t)\|_{H_0^1(\La;w)}^2, 
	\end{split}
\end{equation} 
these bounds imply that all terms in $P_1 z$, $P_2 z$, and $Kz$ belong to $L^2(\La_T)$, so the inner products $(P_1 z, P_2 z)_{L^2(\La_T)}$ and $(P_1 z, Kz)_{L^2(\La_T)}$ are well-defined.

Finally, for any $u\in D(\mathcal A_m^{(2)})$, we have $x_2^\alpha\partial_{x_2}u\in H^1(\Lambda)\hookrightarrow C(\overline\Lambda)$  and $u(0)=u(1)=0$. When we integrate by parts on $(\e,1)\times(0,T)$ and let $\e\to0^+$, all boundary contributions at $x_2=0$ vanish. At  $x_2=1$, the only surviving spatial boundary term in $(P_1z,P_2z)_{L^2(\Lambda_T)}$ is 
\begin{equation*}
	s\int_0^T\Theta(t)|\partial_{x_2}z(1,t)|^2\,\mathrm dt\ge0.
\end{equation*}
All the other spatial boundary terms vanish because they contain $z$ or $z_t$. The nonnegative term above is retained in the identity and may be discarded when taking a lower bound. In particular, noticing that $\be - 1 = \frac{1-\al}{2} > 0$, we have
\begin{equation}\label{11.26.7}
	\lim_{x_2 \to 0^+} (\pt_{x_2}\xi(x_2,t)) \left(x_2^\al \pt_{x_2}z(x_2,t)\right)^2 = 0  \text{ and }  \lim_{x_2 \to 0^+} z_t(x_2,t) \left(x_2^\al \pt_{x_2}z(x_2,t)\right) = 0
\end{equation}
for all $t \in (0,T)$ (respectively, for a.e. $t \in (0,T)$ for the second limit), while all temporal boundary terms at $t = 0$ and $t = T$ vanish due to the exponential decay of $e^{-s\xi}$.

\begin{lemma}\label{11.18.L1}
	Let $u\in D(\mcA_0^{(2)})$, where $\mcA_0^{(2)}$ is defined in \eqref{11.27.1}, and let $q = \f{3\al-1}{2}$ as in \eqref{11.26.1}. Then 
	\begin{equation*}
		\int_\La x_2^q |\pt_{x_2}u|^2\df x_2\leq C\|u\|_{D(\mcA_0^{(2)})}^2,\quad \int_\La x_2^{q-2}|u|^2\df x_2\leq \f{4}{(1-q)^2}\int_\La x_2^q|\pt_{x_2}u|^2\df x_2. 
	\end{equation*}
\end{lemma}

\begin{proof}
Since $u\in D(\mcA_0^{(2)})$, we have $u\in H_0^1(\La;w)$ and $\pt_{x_2}(x_2^\al \pt_{x_2}u)\in L^2(\La)$, which implies $x_2^\al\pt_{x_2}u\in H^1(\La)$. In particular, $x_2^\al \pt_{x_2}u$ is a  continuous function on $\overline{\La}$ with $\|x_2^\al \pt_{x_2}u\|_{L^\iy(\La)}\leq C\|u\|_{D(\mcA_0^{(2)})}$. With $p:=2\al-q = \f{\al+1}{2}$, the condition $p<1$ ensures that $x_2^{-p}$ is integrable on $\La=(0,1)$, and hence
	\begin{equation*}
		\begin{split}
			\int_\La x_2^q|\pt_{x_2}u|^2\df x_2 
			&= \int_\La \left|x_2^{\al} \pt_{x_2}u\right|^2 x_2^{-p}\df x_2 \leq \|x_2^\al \pt_{x_2}u\|_{L^\iy(\La)}^2\int_\La x_2^{-p}\df x_2\\
			&= \f{1}{1-p} \|x_2^\al \pt_{x_2}u\|_{L^\iy(\La)}^2\leq C \|u\|_{D(\mcA_0^{(2)})}^2,
		\end{split}
	\end{equation*}
	where the positive constant $C$ depends only on $\al$ and $q$.
	
	Next, since $\al\in(0,1)$, $u\in D(\mcA_0^{(2)})$, and $u(0)=u(1)=0$, we have for all $x_2\in\La$,
	\begin{equation}\label{11.19.1} 
		|u(x_2)|\leq \int_0^{x_2}s^{-\al}|s^\al \pt_s u(s)|\df s\leq \f{x_2^{1-\al}}{1-\al}\|x_2^\al \pt_{x_2}u\|_{L^\iy(\La)} \leq Cx_2^{1-\al}\|u\|_{D(\mcA_0^{(2)})}.
	\end{equation} 
	Since $q+1-2\al = 1-p > 0$, we obtain
	\begin{equation*}
		0\leq x_2^{q-1}|u(x_2)|^2 \leq C x_2^{q+1-2\al}\|u\|_{D(\mcA_0^{(2)})}^2 \to 0  \quad \text{as } x_2\to 0^+,
	\end{equation*}
	and clearly $\lim_{x_2\to 1^-} x_2^{q-1}|u(x_2)|^2 = 0$. Furthermore, on $\La=(0,1)$ we have
	\begin{equation*}
		x_2^{q-2}|u(x_2)|^2\leq Cx_2^{q-2\al} \|u\|_{D(\mcA_0^{(2)})}^2=Cx_2^{-p}\|u\|_{D(\mcA_0^{(2)})}^2\in L^1(\La),
	\end{equation*}
	which guarantees that $\int_\La x_2^{q-2}|u|^2\df x_2 < +\iy$. Integrating by parts and applying the Cauchy-Schwarz inequality, we deduce that
	\begin{equation*}
		\begin{split}
			\int_\La x_2^{q-2}|u|^2\df x_2 
			&= \f{1}{q-1}\int_\La \pt_{x_2}(x_2^{q-1}) |u|^2\df x_2 
			= -\f{2}{q-1}\int_\La x_2^{q-1} {u}\,\pt_{x_2}u\df x_2 \\
			&\leq \f{2}{|1-q|} \left(\int_\La x_2^{q-2}|u|^2\df x_2\right)^{\f{1}{2}} \left(\int_\La x_2^q |\pt_{x_2}u|^2\df x_2\right)^{\f{1}{2}},
		\end{split}
	\end{equation*}
	which gives
	\begin{equation*}
		\int_\La x_2^{q-2}|u|^2\df x_2 \leq \f{4}{(1-q)^2}\int_\La x_2^q|\pt_{x_2}u|^2\df x_2.
	\end{equation*}
	This completes the proof. 
\end{proof}

\subsection{Computation of the inner products}\label{S3.2}

Based on the justification in Section \ref{S3.1} and the boundary behavior (see \eqref{11.19.1})
\begin{equation*}
	0 \leq \lim_{x_2 \to 0^+} x_2^{\f{3(\al-1)}{2}} z^2(x_2,t) \leq C \lim_{x_2 \to 0^+} x_2^{\f{1-\al}{2}} \|z(\cdot,t)\|_{D(\mcA_0^{(2)})}^2 = 0,
\end{equation*}
a direct integration by parts yields
\begin{equation*}
	\begin{split}
		(P_1z,P_2z)_{L^2(\La_T)}
		&= s\int_0^T \left.(\pt_{x_2}\xi)(x_2^\al \pt_{x_2}z)^2\right|_{x_2=0}^{x_2=1}\df t + \f{s}{2}\iint_{\La_T} z^2\pt_{x_2}\left[x_2^\al \pt_{x_2x_2}(x_2^\al \pt_{x_2}\xi)\right]\df x_2\,\df t\\
		&\quad - 2s\iint_{\La_T} x_2^\al (\pt_{x_2}z)^2\pt_{x_2}(x_2^\al \pt_{x_2}\xi)\df x_2\,\df t + \al s\iint_{\La_T} (\pt_{x_2}z)^2x_2^{2\al-1} \pt_{x_2}\xi\df x_2\,\df t\\
		&\quad - \f{s}{2}\iint_{\La_T} z^2\xi_{tt}\df x_2\,\df t + 2s^2\iint_{\La_T} z^2x_2^\al (\pt_{x_2}\xi)(\pt_{x_2}\xi_t)\df x_2\,\df t\\
		&\quad - s^3\iint_{\La_T} z^2 x_2^\al (\pt_{x_2}\xi)\left(\pt_{x_2}\left[x_2^\al (\pt_{x_2}\xi)^2\right]\right)\df x_2\,\df t, 
	\end{split}
\end{equation*}
and 
\begin{equation*}
	(P_1z,Kz)_{L^2(\La_T)}
	= \al\mu_ms\iint_{\La_T}z^2x_2^{2\al-1}\pt_{x_2}\xi\df x_2\,\df t. 
\end{equation*}

\subsection{Estimates}\label{S3.3}

(1) Direct calculations give
\begin{equation*}
	|\Theta'|\leq 4T\Theta^\f{5}{4},\quad |\Theta\Theta'|\leq 4T\Theta^\f{9}{4},\quad |\Theta''|\leq 20T^2\Theta^\f{3}{2}. 
\end{equation*}

(2) Set
\begin{equation*}
	\La^1=(0,c), \quad \La^2=(c,d),\quad \La^3=(d,1), \quad \La_T^i=\La^i\ts (0,T)  \text{ for } i=1,2,3. 
\end{equation*}

(I) On $\La_T^1$: noting from \eqref{11.26.1} that $\be-1=\f{1-\al}{2}$, $\al+\be-1=h$, $2\al+\be-2=q$, $\al+\be-2=h-1$, and $2\al+\be-3=\f{3(\al-1)}{2}$, we compute
\begin{equation*}
	\pt_{x_2} \eta=\be x_2^{\be-1}, \quad \pt_{x_2} \xi=-\be\Theta x_2^{\be-1}, \quad x_2^\al \pt_{x_2}\xi=-\be \Theta x_2^{h},\quad  \pt_{x_2}\xi_t=-\be\Theta'x_2^{\be-1}, 
\end{equation*}
\begin{equation*}
	x_2^{2\al-1}\pt_{x_2}\xi=-\be \Theta x_2^{q}, \quad \pt_{x_2}(x_2^\al \pt_{x_2}\xi)=-\be h\Theta x_2^{h-1}, \quad \xi_{tt}=\Theta''(\ga-\eta), 
\end{equation*}
\begin{equation*}
	x_2^\al \pt_{x_2x_2}(x_2^\al\pt_{x_2}\xi)=-\be h(h-1)\Theta x_2^\f{3(\al-1)}{2}, \quad x_2^\al (\pt_{x_2}\xi)(\pt_{x_2}\xi_t)=\be^2\Theta\Theta' x_2, 
\end{equation*}
and  
\begin{equation*}
	\pt_{x_2}\left[x_2^\al \pt_{x_2x_2}(x_2^\al\pt_{x_2}\xi)\right]=-\be h(h-1)\f{3(\al-1)}{2}\Theta x_2^{q-2}, \ x_2^\al (\pt_{x_2}\xi)\left(\pt_{x_2}\left[x_2^\al (\pt_{x_2}\xi)^2\right]\right)=-\be^3\Theta^3 x_2^{h}. 
\end{equation*}

(II) On $\La_T^3$: since $\zeta(x_2)=0$, we have $\eta(x_2)=(1-x_2)x_2^{-\al}$ and
\begin{equation*}
	\pt_{x_2}\eta(x_2) = -x_2^{-\al} - \al(1-x_2)x_2^{-\al-1} = -\left[(1-\al)x_2^{-\al} + \al x_2^{-\al-1}\right] < 0,
\end{equation*}
which implies
\begin{equation*}
	\pt_{x_2}\xi(x_2,t) = \Theta(t)\left[(1-\al)x_2^{-\al} + \al x_2^{-\al-1}\right] > 0.
\end{equation*}
Moreover,
\begin{equation*}
	x_2^\al \pt_{x_2}\xi = \Theta \left(1-\al + \al x_2^{-1}\right), \quad \pt_{x_2}\xi_t = \Theta'\left[(1-\al)x_2^{-\al} + \al x_2^{-\al-1}\right], 
\end{equation*}
\begin{equation*}
	\xi_{tt} = \Theta''(\ga-\eta), \quad \pt_{x_2}\left(x_2^\al \pt_{x_2}\xi\right) = -\al \Theta x_2^{-2},\quad x_2^\al \pt_{x_2x_2}\left(x_2^\al \pt_{x_2}\xi\right) = 2\al \Theta x_2^{\al-3},
\end{equation*}
\begin{equation*}
	x_2^{2\al-1}\pt_{x_2}\xi = \Theta\left[(1-\al)x_2^{\al-1} + \al x_2^{\al-2}\right], \quad \pt_{x_2}\left[x_2^\al \pt_{x_2x_2} (x_2^\al \pt_{x_2}\xi)\right]=-2\al (3-\al)\Theta x_2^{\al-4},
\end{equation*}
\begin{equation*}
	x_2^\al (\pt_{x_2}\xi)(\pt_{x_2}\xi_t) = \Theta\Theta'\left[(1-\al)^2 x_2^{-\al} + 2\al(1-\al)x_2^{-\al-1} + \al^2 x_2^{-\al-2}\right],
\end{equation*}
and
\begin{equation*}
	x_2^\al (\pt_{x_2}\xi)\left(\pt_{x_2}\left[x_2^\al (\pt_{x_2}\xi)^2\right]\right) = -\al \Theta^3 x_2^{-\al-1}\left[(1-\al)+\al x_2^{-1}\right]^2\left[(1-\al)+(2+\al)x_2^{-1}\right] < 0.
\end{equation*}

(3) By \eqref{11.26.2} and \eqref{11.26.7}, the boundary term in $(P_1z,P_2z)_{L^2(\La_T)}$ simplifies to
\begin{equation}\label{11.26.9}
	s\int_0^T \left.(\pt_{x_2}\xi)(x_2^\al \pt_{x_2}z)^2\right|_{x_2=0}^{x_2=1}\df t = s\int_0^T \Theta(t) |\pt_{x_2}z(1,t)|^2\df t.
\end{equation}

(4) We decompose the volume integrals in $(P_1z,P_2z)_{L^2(\La_T)}$ into three parts corresponding to $\La_T^1$, $\La_T^2$, and $\La_T^3$:
\begin{equation*}
	(P_1z,P_2z)_{L^2(\La_T)} = s\int_0^T \Theta(t) |\pt_{x_2}z(1,t)|^2\df t + I_1 + I_2 + I_3.
\end{equation*}

(5) By (I), the integral $I_1$ on $\La_T^1$ is given by
\begin{equation*}
	\begin{split}
		I_1
		&=\be s\iint_{\La_T^1}\Theta x_2^q(\pt_{x_2}z)^2\df x_2\df t+\be^3s^3\iint_{\La_T^1}\Theta^3 x_2^h z^2\df x_2\df t\\
		&\hspace{4.5mm}-\f{s}{2}\iint_{\La_T^1}\Theta'' z^2(\ga-\eta)\df x_2\df t+2\be^2s^2\iint_{\La_T^1}\Theta\Theta' z^2 x_2\df x_2\df t\\
		&\hspace{4.5mm}-\f{\be}{2} h(h-1)\f{3(\al-1)}{2}s\iint_{\La_T^1} \Theta z^2 x_2^{q-2}\df x_2\df t.
	\end{split}
\end{equation*}
Note that on $\La_T^1$, since $0 < x_2 < c < 1$ and $\f{(q-2)+h}{2} = \al - 1 < 0$, we have $x_2^{\al-1} > 1$. Applying Cauchy's inequality, Lemma \ref{11.18.L1}, and the bound $1 \leq \left(\f{T}{2}\right)^8 \Theta$ yields
\begin{equation*}
	\begin{split}
		\f{s}{2}\iint_{\La_T^1} |\Theta''|z^2(\ga-\eta)\df x_2\,\df t
		&\leq 10T^2 \ga s\iint_{\La_T^1}\Theta^{\f{3}{2}} z^2 x_2^{\al-1}\df x_2\,\df t\\
		&\leq 5T^6\ga^{\f{2}{3}} \iint_{\La_T^1} \Theta z^2 x_2^{q-2} \df x_2\,\df t + 5T^{-2}\ga^{\f{4}{3}}s^2\iint_{\La_T^1} \Theta^2z^2 x_2^h\df x_2\,\df t\\
		&\leq C T^6 \iint_{\La_T}\Theta x_2^q(\pt_{x_2}z)^2\df x_2\,\df t + C T^6 s^2\iint_{\La_T^1} \Theta^3z^2x_2^h\df x_2\,\df t,
	\end{split}
\end{equation*}
where $C > 0$ depends only on $\al$, $c$, and $d$. Similarly, since $x_2 \leq x_2^h$ on $\La_T^1$ and $|\Theta\Theta'| \leq 4T\Theta^{\f{9}{4}}$, we obtain
\begin{equation*}
	2\be^2 s^2 \iint_{\La_T^1} |\Theta\Theta'| z^2 x_2\,\df x_2\,\df t 
	\leq C T^7 s^2 \iint_{\La_T^1} \Theta^3 z^2 x_2^h\,\df x_2\,\df t.
\end{equation*}
Moreover, since $q < 1$ (so the boundary term at $x_2 = c$ in integration by parts has a favorable non-positive sign), the same argument as in Lemma \ref{11.18.L1} yields
\begin{equation}\label{11.27.4}
	\int_0^c x_2^{q-2}z^2\,\df x_2 \leq \f{4}{(1-q)^2}\int_0^c x_2^q(\pt_{x_2}z)^2\,\df x_2.
\end{equation}
Since $\f{1}{2}\left|\be h(h-1)\f{3(\al-1)}{2}\right| \cdot \f{4}{(1-q)^2} = \f{1+\al}{3}\be$, this directly gives
\begin{equation*}
	\f{1}{2}\left|\be h(h-1)\f{3(\al-1)}{2}\right| s\iint_{\La_T^1} \Theta z^2 x_2^{q-2}\,\df x_2\,\df t 
	\leq \f{(1+\al)\be s}{3}\iint_{\La_T^1} \Theta x_2^q(\pt_{x_2}z)^2\,\df x_2\,\df t.
\end{equation*}
Combining the above estimates (with the $\Theta''$-term also restricted to $\La_T^1$ via the same Hardy inequality on $(0,c)$) and using $T^7 \leq \f{1}{2}(T^6 + T^8)$, we conclude that there exists a constant $\la_1 > 0$ (depending only on $\al$, $c$, and $d$) such that for all $s \geq \la_1 (T^6 + T^8)$,
\begin{equation}\label{11.27.2}
	I_1 \geq \f{\be(2-\al)}{6}s \iint_{\La_T^1} \Theta x_2^q (\pt_{x_2}z)^2\,\df x_2\,\df t + \f{\be^3 s^3}{2} \iint_{\La_T^1} \Theta^3 x_2^h z^2\,\df x_2\,\df t.
\end{equation}

(6) Similarly, by (II), the integral $I_3$ on $\La_T^3$ can be expressed as
\begin{equation*}
	\begin{split}
		I_3
		&=\alpha s\iint_{\Lambda_T^3}\Theta(\partial_{x_2}z)^2[(1-\alpha)x_2^{\alpha-1}+(2+\alpha)x_2^{\alpha-2}]\,\mathrm dx_2\,\mathrm dt\\
		&\quad+\alpha s^3\iint_{\Lambda_T^3}\Theta^3x_2^{-\alpha-1}[(1-\alpha)+\alpha x_2^{-1}]^2[(1-\alpha)+(2+\alpha)x_2^{-1}]z^2\,\mathrm dx_2\,\mathrm dt\\
		&\quad-\alpha(3-\alpha)s\iint_{\Lambda_T^3}\Theta x_2^{\alpha-4}z^2\,\mathrm dx_2\,\mathrm dt-\frac{s}{2}\iint_{\Lambda_T^3}\Theta''(\gamma-\eta)z^2\,\mathrm dx_2\,\mathrm dt\\
		&\quad+2s^2\iint_{\Lambda_T^3}\Theta\Theta'z^2[(1-\alpha)^2x_2^{-\alpha}+2\alpha(1-\alpha)x_2^{-\alpha-1}+\alpha^2x_2^{-\alpha-2}]\,\mathrm dx_2\,\mathrm dt.
	\end{split}
\end{equation*}
Since $\al \in (0,1)$ and $d < x_2 < 1$ on $\La_T^3$ (so that $x_2^{-1} > 1$), we observe that
\begin{equation*}
	\al\left[(1-\al)x_2^{\al-1} + (2+\al)x_2^{\al-2}\right] \geq 3\al
\end{equation*}
and
\begin{equation*}
	\al x_2^{-\al-1}\left[(1-\al)+\al x_2^{-1}\right]^2\left[(1-\al)+(2+\al)x_2^{-1}\right] \geq 3\al.
\end{equation*}
Using the estimates for $\Theta'$ and $\Theta''$ together with $d < x_2 < 1$, we obtain
\begin{equation*}
	\begin{split}
		&\left| \al(3-\al) s \iint_{\La_T^3} \Theta z^2 x_2^{\al-4} \df x_2\,\df t \right| + \left| \f{s}{2}\iint_{\La_T^3} \Theta'' z^2(\ga-\eta)\df x_2\,\df t \right| \\
		&\quad + \left| 2s^2 \iint_{\La_T^3} \Theta\Theta' z^2 \left[ (1-\al)^2 x_2^{-\al} + 2\al(1-\al)x_2^{-\al-1} + \al^2 x_2^{-\al-2} \right] \df x_2\,\df t \right| \\
		&\leq C \iint_{\La_T^3} \left( s\Theta + s T^2 \Theta^{\f{3}{2}} + s^2 T \Theta^{\f{9}{4}} \right) z^2 \df x_2\,\df t,
	\end{split}
\end{equation*}
where $C > 0$ depends only on $\al$ and $d$. Noting that $1 \leq \left(\f{T}{2}\right)^8 \Theta$ and $T^7 \leq \f{1}{2}(T^6 + T^8)$, we conclude that there exists a constant $\la_3 > 0$ (depending only on $\al$ and $d$) such that for all $s \geq \la_3 (T^6 + T^8)$, we have
\begin{equation}\label{11.27.3}
	I_3 \geq 3\al s \iint_{\La_T^3} \Theta (\pt_{x_2}z)^2 \df x_2\,\df t + \al s^3 \iint_{\La_T^3} \Theta^3 z^2 \df x_2\,\df t.
\end{equation}

(7) From (4), \eqref{11.27.2}, and \eqref{11.27.3}, there exists a constant $\la > 0$ (depending only on $\al$, $c$, and $d$) such that for all $s \geq \la(T^6 + T^8)$, we obtain
\begin{equation*}
	\begin{split} 
	(P_1z,P_2z)_{L^2(\La_T)}
	&\geq C\iint_{\La_T} \left(s\Theta x_2^q(\pt_{x_2}z)^2 + s^3\Theta^3 x_2^h z^2 \right)\df x_2\,\df t \\
	&\hspace{4.5mm}- C\iint_{\La_T^2} \left( s\Theta (\pt_{x_2}z)^2 + s^3\Theta^3 z^2 \right) \df x_2\,\df t,
	\end{split} 
\end{equation*}
where $C > 0$ depends only on $\al$, $c$, and $d$.

(8) Fix the cutoff in the definition of $\eta$ by rescaling a fixed smooth cutoff profile to $(c,d)$, and set
\begin{equation*}
	M_\eta:=\sup_{0<x_2\le1}x_2^{1-\beta}|\eta'(x_2)|<\infty,\quad H_q:=\frac{4}{(1-q)^2},\quad \vartheta:=\frac{\beta-1}{\beta}\in(0,1).
\end{equation*} 
Here $M_\eta$ depends only on $\alpha,c,d$. Since $q=\vartheta(q-2)+(1-\vartheta)h$, H\"older's inequality and Lemma~\ref{11.18.L1} imply, with
\begin{equation*}
	G_0(t):=\int_\Lambda x_2^q|\partial_{x_2}z|^2\,\mathrm dx_2,\quad H_0(t):=\int_\Lambda x_2^h z^2\,\mathrm dx_2,
\end{equation*}
that
\begin{equation*}
	\int_\Lambda x_2^qz^2\,\mathrm dx_2\le\left(\int_\Lambda x_2^{q-2}z^2\,\mathrm dx_2\right)^\vartheta H_0(t)^{1-\vartheta}\le H_q^\vartheta G_0(t)^\vartheta H_0(t)^{1-\vartheta}.
\end{equation*}
Consequently, for every $\e>0$ and every $\mu_m\ge0$, Young's inequality gives
\begin{equation*}
	\mu_m\int_\Lambda x_2^qz^2\,\mathrm dx_2\le\e  G_0(t)+\frac{\vartheta^{\beta-1}}{\beta}H_q^{\beta-1}\e^{-(\beta-1)}\mu_m^\beta H_0(t).
\end{equation*}
The exact potential cross term is
\begin{equation*}
	(P_1z,Kz)_{L^2(\Lambda_T)}=-\alpha\mu_ms\iint_{\Lambda_T}\Theta x_2^{2\alpha-1}\eta'(x_2)z^2\,\mathrm dx_2\,\mathrm dt\ge-\alpha M_\eta\mu_ms\iint_{\Lambda_T}\Theta x_2^qz^2\,\mathrm dx_2\,\mathrm dt.
\end{equation*}
Using the preceding interpolation estimate, rescaling $\e$, and using $\Theta\le 2^{-16}T^{16}\Theta^3$, we obtain
\begin{equation*}
	\begin{split}
		(P_1z,Kz)_{L^2(\Lambda_T)}
		&\ge-\e s\iint_{\Lambda_T}\Theta x_2^q|\partial_{x_2}z|^2\,\mathrm dx_2\,\mathrm dt-C_\e\mu_m^\beta T^{16}s\iint_{\Lambda_T}\Theta^3x_2^hz^2\,\mathrm dx_2\,\mathrm dt,
	\end{split}
\end{equation*}
where $C_\e$ depends only on $\e,\alpha,c,d$. By taking $\e$ sufficiently small and then $\lambda_0$ sufficiently large, we can absorb both negative terms into the positive terms in (7) whenever $s=\lambda[T^6+T^8(1+\mu_m^\f{\be}{2})]$ with $\lambda\ge\lambda_0$. Indeed, $s^2\ge\lambda^2T^{16}\mu_m^\beta$.

(9) Taking $s=\la (T^6+T^8(1+\mu_m^{\f{\be}{2}}))$ with $\la\geq \la_0$ (for some constant $\la_0\geq 1$ depending only on $\al, c$ and $d$), we deduce from (7) and (8) that
\begin{equation*}
	\begin{split}
		&s\iint_{\La_T} \Theta x_2^q(\pt_{x_2}z)^2\df x_2\df t+s^3\iint_{\La_T} \Theta^3 x_2^h z^2\df x_2\df t\\
		&\leq C\|e^{-s\xi}f_m\|_{L^2(\La_T)}^2+C\iint_{\La_T^2} \left(s\Theta (\pt_{x_2}z)^2+s^3\Theta^3z^2\right)\df x_2\df t, 
	\end{split}
\end{equation*}
where the positive constants $C$ depend only on $\al, c$ and $d$.

(10) Since 
\begin{equation*}
	\pt_{x_2} \vp_m=\pt_{x_2}\left(ze^{s\xi}\right)=e^{s\xi}\left(\pt_{x_2}z+sz\pt_{x_2}\xi\right), 
\end{equation*}
we have 
\begin{equation*}
	e^{-2s\xi}x_2^q|\partial_{x_2}\varphi_m|^2\le 2x_2^q|\partial_{x_2}z|^2+2M_\eta^2s^2\Theta^2x_2^hz^2.
\end{equation*}
Hence, from (9) and $\pt_{x_2}z=e^{-s\xi}(\pt_{x_2}\vp_m-s\vp_m\pt_{x_2}\xi)$, 
we obtain for $s=\la (T^6+T^8(1+\mu_m^{\f{\be}{2}}))$ with $\la\geq \la_0$, 
\begin{equation*}
	\begin{split}
		&s\iint_{\La_T} e^{-2s\xi}\Theta x_2^q(\pt_{x_2}\vp_m)^2\df x_2\df t+s^3\iint_{\La_T} e^{-2s\xi}\Theta^3 x_2^h \vp_m^2\df x_2\df t\\
		&\leq C\|e^{-s\xi}f_m\|_{L^2(\La_T)}^2+C\iint_{\La_T^2} e^{-2s\xi}\left(s\Theta (\pt_{x_2}\vp_m)^2+s^3\Theta^3\vp_m^2\right)\df x_2\df t, 
	\end{split}
\end{equation*}
where the positive constants $C$ depend only on $\al, c$ and $d$.

(11) From Lemma \ref{11.18.L1}, for all $u\in D(\mcA_0^{(2)})$, 
\begin{equation*}
	\begin{split}
		\int_\La u^2\df x_2
		&\leq \int_\La x_2^{\al-1}u^2\df x_2\leq \left(\int_\La x_2^hu^2\df x_2\right)^{\f{1}{2}}\left(\int_\La x_2^{q-2}u^2\df x_2\right)^{\f{1}{2}}\\
		&\leq \f{2}{1-q}\left(\int_\La x_2^hu^2\df x_2\right)^{\f{1}{2}}\left(\int_\La x_2^q(\pt_{x_2}u)^2\df x_2\right)^{\f{1}{2}}. 
		\end{split}
	\end{equation*} 
Hence, taking $u=z=\vp_m e^{-s\xi}$, we obtain
\begin{equation*}
	\begin{split}
		&s^2\iint_{\La_T} \Theta^2\vp_m^2e^{-2s\xi}\df x_2\df t\leq \f{1}{1-q}s\iint_{\La_T}\Theta  x_2^q(\pt_{x_2}z)^2\df x_2\df t+\f{1}{1-q}s^3\iint_{\La_T} \Theta^3 x_2^h z^2\df x_2\df t. 
	\end{split}
\end{equation*}
Combining this with (9) and the preceding estimate yields 
\begin{equation*}
	\begin{split}
		&s\iint_{\La_T} \Theta x_2^q(\pt_{x_2}\vp_m)^2e^{-2s\xi}\df x_2\df t+s^2\iint_{\La_T} \Theta^2\vp_m^2e^{-2s\xi}\df x_2\df t +s^3\iint_{\La_T} e^{-2s\xi}\Theta^3 x_2^h \vp_m^2\df x_2\df t\\
		&\leq C\|e^{-s\xi}f_m\|_{L^2(\La_T)}^2+C\iint_{\La_T^2} e^{-2s\xi}\left(s\Theta (\pt_{x_2}\vp_m)^2+s^3\Theta^3\vp_m^2\right)\df x_2\df t.
	\end{split}
\end{equation*}

(12) To eliminate the local gradient term on $\La_T^2=(c,d)\ts (0,T)$, we pick a cutoff function $\zeta \in C_0^\infty(a,b)$ satisfying $0 \leq \zeta \leq 1$ and $\zeta \equiv 1$ on $(c,d)$. Note that $x_2^\al \geq a^\al > 0$ on $(a,b)$. Multiplying the equation \eqref{11.15.21} satisfied by $\vp_m$ by $\zeta^2 s\Theta e^{-2s\xi}\vp_m$, integrating over $\La_T$, and using integration by parts together with $\lim_{t\to 0^+,T^-}(\Theta e^{-2s\xi})(x_2,t)=0$ and $\mu_m \geq 0$, we obtain
\begin{equation*}
	\begin{split}
		&\f{a^\al}{2}s\iint_{\La_T^2} \Theta(\pt_{x_2}\vp_m)^2e^{-2s\xi}\df x_2\df t\\
		&\leq \f{1}{2}s\iint_{\La_T} \Theta\zeta^2  x_2^\al  (\pt_{x_2}\vp_m)^2e^{-2s\xi}\df x_2\df t\\
		&\leq C\|e^{-s\xi}f_m\|_{L^2(\La_T)}^2 + C\iint_{(a,b)\ts (0,T)} \left(s\Theta + s^2\Theta^2 + s^3\Theta^3\right)\vp_m^2e^{-2s\xi}\df x_2\df t\\
		&\leq C\|e^{-s\xi}f_m\|_{L^2(\La_T)}^2 + C\iint_{(a,b)\ts (0,T)} s^3\Theta^3\vp_m^2e^{-2s\xi}\df x_2\df t,
	\end{split}
\end{equation*}
where we used $s\Theta + s^2\Theta^2 \leq C s^3\Theta^3$ for $s=\la (T^6+T^8(1+\mu_m^{\f{\be}{2}}))$ with $\la \geq \la_0$. Consequently, 
\begin{equation*}
	\begin{split} 
		&\iint_{\La_T^2} e^{-2s\xi}\left(s\Theta (\pt_{x_2}\vp_m)^2+s^3\Theta^3\vp_m^2\right)\df x_2\df t \\
		&\leq C\|e^{-s\xi}f_m\|_{L^2(\La_T)}^2 + C\iint_{(a,b)\ts (0,T)} e^{-2s\xi}s^3\Theta^3\vp_m^2\df x_2\df t, 
	\end{split} 
\end{equation*}
where the positive constants $C$ depend only on $\al, a, b, c$ and $d$. Combining this with (11), we obtain for $s=\la(T^6+T^8(1+\mu_m^{\f{\be}{2}}))$ with $\la\geq \la_0$,
\begin{equation*}
	\begin{split}
		&s\iint_{\La_T} \Theta x_2^q(\pt_{x_2}\vp_m)^2e^{-2s\xi}\df x_2\df t+s^2\iint_{\La_T} \Theta^2\vp_m^2e^{-2s\xi}\df x_2\df t+s^3\iint_{\La_T} \Theta^3 x_2^h \vp_m^2e^{-2s\xi}\df x_2\df t\\
		&\leq C\|e^{-s\xi}f_m\|_{L^2(\La_T)}^2+Cs^3\iint_{(a,b)\ts (0,T)}\Theta^3\vp_m^2e^{-2s\xi}\df x_2\df t,
	\end{split}
\end{equation*}
where the positive constants $C$ depend only on $\al, a, b, c$ and $d$. This proves the following theorem. 

\begin{theorem}\label{11.27.T1}
	Let $\al\in (0,1)$, $\vp_m^0\in D(\mcA_0^{(2)})$, and $f_m\in H^1(0,T; H_0^1(\La;w))$, and let $\vp_m$ be the corresponding solution to \eqref{11.15.21} with $\mu_m\geq 0$. Suppose that $0<a<c<d<b<1$, and let $\Theta$ and $\xi$ be defined in \eqref{11.26.3}. Then there exist positive constants $\la_0\geq 1$ and $C>0$, depending only on $\al, a, b, c$, and $d$, such that for all $\la\geq \la_0$ and $s=\la(T^6+T^8(1+\mu_m^{\f{\be}{2}}))$,
	\begin{equation}\label{11.27.5}
		\begin{split}
			&s\iint_{\La_T} \Theta x_2^q(\pt_{x_2}\vp_m)^2e^{-2s\xi}\df x_2\df t+s^3\iint_{\La_T} \Theta^3 x_2^h \vp_m^2e^{-2s\xi}\df x_2\df t\\
			&\hspace{4.5mm}+s^2\iint_{\La_T} \Theta^2\vp_m^2e^{-2s\xi}\df x_2\df t \\
			&\leq C\|e^{-s\xi}f_m\|_{L^2(\La_T)}^2+Cs^3\iint_{(a,b)\ts (0,T)}\Theta^3\vp_m^2e^{-2s\xi}\df x_2\df t.
		\end{split}
	\end{equation}
\end{theorem}

\begin{corollary}\label{11.27.C1}
	Let $\al\in (0,1)$, $\vp_m^0\in L^2(\La)$, and $f_m=0$, and let $\vp_m$ be the corresponding mild solution to \eqref{11.15.21} with $\mu_m\geq 0$. Suppose that $0<a<c<d<b<1$. Then there exists a positive constant $C$, depending only on $\al, a, b, c$, and $d$, such that  
	\begin{equation}\label{11.27.6}
		\int_\La \vp_m^2(T)\df x_2\leq Ce^{C\left(1+T^{-2}+\mu_m^{\f{\be}{2}}\right)}\iint_{(a,b)\ts (0,T)} \vp_m^2\df x_2\df t. 
	\end{equation}
\end{corollary}

\begin{proof}[Proof of Corollary \ref{11.27.C1}]
	By a standard density argument (since $D(\mcA_0^{(2)})$ is dense in $L^2(\La)$ and the mild solution depends continuously on the initial data in $C([0,T];L^2(\La))$), it suffices to prove \eqref{11.27.6} for $\vp_m^0\in D(\mcA_0^{(2)})$.
	
	Multiplying \eqref{11.15.21} (with $f_m=0$) by $\vp_m$, integrating over $\La$, and using integration by parts together with $\mu_m\geq 0$, we obtain
	\begin{equation*}
		\frac12\frac{\mathrm d}{\mathrm dt}\int_\Lambda\varphi_m^2(x_2,t)\,\mathrm dx_2+\int_\Lambda x_2^\alpha|\partial_{x_2}\varphi_m|^2\,\mathrm dx_2+\mu_m\int_\Lambda x_2^\alpha\varphi_m^2\,\mathrm dx_2=0.
	\end{equation*}
	This implies that $t \mapsto \int_\La \vp_m^2(x_2,t)\df x_2$ is nonincreasing on $[0,T]$. In particular, for all $t\in [\f{T}{4},\f{3T}{4}]$,
	\begin{equation}\label{11.27.7}
		\int_\La \vp_m^2(x_2,T)\df x_2 \leq \int_\La \vp_m^2(x_2,t)\df x_2.
	\end{equation}
	
	Next, recall from \eqref{11.26.3} that $\Theta(t)=[t(T-t)]^{-4}$ and $\xi(x_2,t)=\Theta(t)(\ga-\eta(x_2))$ with $\ga=\|\eta\|_{L^\iy(\La)}+1$. Note that $1\leq \ga-\eta(x_2)\leq 2\|\eta\|_{L^\iy(\La)}+1 =: C_0$ for all $x_2\in \La$, and for all $t\in [\f{T}{4},\f{3T}{4}]$,
	\begin{equation*}
		\f{2^8}{T^8} \leq \Theta(t) \leq \f{2^{16}}{3^4 T^8}.
	\end{equation*}
	Fixing $\la=\la_0$ (where $\la_0\geq 1$ is given by Theorem \ref{11.27.T1}) and taking $s=\la_0(T^6+T^8(1+\mu_m^{\f{\be}{2}}))$, we have $s\geq \la_0 T^8$. Consequently, on $\La\ts[\f{T}{4},\f{3T}{4}]$,
	\begin{equation*}
		s^2\Theta^2(t)e^{-2s\xi(x_2,t)} \geq (\la_0 T^8)^2\left(\f{2^8}{T^8}\right)^2 e^{-2C_0\la_0(T^6+T^8(1+\mu_m^{\f{\be}{2}}))\f{2^{16}}{3^4 T^8}} \geq c_1 e^{-C\left(1+T^{-2}+\mu_m^{\f{\be}{2}}\right)}
	\end{equation*}
	for some positive constants $c_1$ and $C$ depending only on $\al, a, b, c$, and $d$.
	
Integrating \eqref{11.27.7} over $t\in [\f{T}{4},\f{3T}{4}]$ and combining the result with the above lower bound, we deduce that
	\begin{equation*}
		\begin{split}
			\int_\La \vp_m^2(x_2,T)\df x_2 
			&\leq \f{2}{T}\int_{\f{T}{4}}^{\f{3T}{4}}\int_\La \vp_m^2(x_2,t)\df x_2\df t \leq \f{C}{T} e^{C\left(1+T^{-2}+\mu_m^{\f{\be}{2}}\right)} s^2\iint_{\La_T} \Theta^2\vp_m^2 e^{-2s\xi}\df x_2\df t.
		\end{split}
	\end{equation*}
	Applying Theorem \ref{11.27.T1} with $f_m=0$ yields
	\begin{equation*}
		\int_\La \vp_m^2(x_2,T)\df x_2 \leq \f{C}{T} e^{C\left(1+T^{-2}+\mu_m^{\f{\be}{2}}\right)} s^3\iint_{(a,b)\ts (0,T)} \Theta^3\vp_m^2 e^{-2s\xi}\df x_2\df t.
	\end{equation*}
	Finally, using $\xi(x_2,t)\geq \Theta(t)>0$ on $\La_T$ together with $\sup_{y>0} y^3 e^{-2y} < \iy$ and $T^{-1}\leq C e^{C T^{-2}}$ for all $T>0$, we conclude that
	\begin{equation*}
		\int_\La \vp_m^2(x_2,T)\df x_2 \leq C e^{C\left(1+T^{-2}+\mu_m^{\f{\be}{2}}\right)}\iint_{(a,b)\ts (0,T)} \vp_m^2\df x_2\df t,
	\end{equation*}
	which completes the proof.
\end{proof}

\subsection{Controllability}\label{S3.4}

\begin{lemma}[\cite{Rousseau}]\label{11.27.L1}
	Let $0<a_1<b_1<1$. Then there exist positive constants $C_1$ and $C_2$, depending only on $a_1$ and $b_1$, such that for all $n\in\N^*$ and all $c_1,\dots, c_n\in\R$, 
	\begin{equation}\label{11.27.8}
		\sum_{k=1}^n |c_k|^2\leq C_1e^{C_2n}\int_{(a_1,b_1)}\left|\sum_{k=1}^n c_k\sqrt{2}\sin (k\pi x_1)\right|^2\df x_1. 
	\end{equation}
	For example, we can take
	\begin{equation*}
		C_1=\f{2}{b_1-a_1},\quad C_2=\log 3+4\log\left(\f{8e}{b_1-a_1}\right). 
	\end{equation*}
\end{lemma}

Define
\begin{equation*}
	\mcE_n=\left\{\sum_{k=1}^ng_k(x_2)X_k(x_1)\colon g_1,\dots,g_n\in L^2(\La)\right\},
\end{equation*}
where $X_k(x_1)=\sqrt{2}\sin(k\pi x_1)$. Note that $\mcE_n$ is an infinite-dimensional closed subspace of $L^2(\Om)$. Let 
\begin{equation*}
	P_n: L^2(\Om)\ra \mcE_n, \quad P_n u(x_1,x_2)=\sum_{k=1}^n \left(\int_\La  u(y_1,x_2)X_k(y_1)\df y_1\right)X_k(x_1)
\end{equation*}
denote the orthogonal projection from $L^2(\Om)$ onto $\mcE_n$.

\begin{lemma}\label{11.27.L2}
	Let $T>0$, $\al\in (0,1)$, and $\om=(a_1,b_1)\ts (a_2,b_2)$ with $0<a_1<b_1<1$ and $0<a_2<b_2<1$. Then there exist positive constants $C_1, C_2$, and $C_3$, depending only on $\al, \om$, and $T$, such that for all $0\leq t_1<t_2\leq T$ with $\tau:=t_2-t_1$, all $n\in\N^*$, and all $\vp^1\in L^2(\Om)$, there exists a control $f_n\in L^2(\om\ts (t_1,t_2))$ for which the corresponding solution $\vp$ to \eqref{11.15.1} on $(t_1,t_2)$ with initial datum $\vp(t_1)=\vp^1$ satisfies
	\begin{equation}\label{11.27.12}
		P_n\vp(t_2)=0,
	\end{equation}
	and 
	\begin{equation}\label{11.27.10}
		\|f_n\|_{L^2(\om\ts (t_1,t_2))}^2\leq C_1e^{C_2\tau^{-2}+C_3n^\be}\|\vp^1\|_{L^2(\Om)}^2,  
	\end{equation}
	and 
	\begin{equation}\label{11.27.11}
		\|\vp(t_2)\|_{L^2(\Om)}^2\leq 2\left(1+C_1Te^{C_2\tau^{-2}+C_3n^\be}\right)\|\vp^1\|_{L^2(\Om)}^2.
	\end{equation}
\end{lemma}

\begin{proof}
	We prove this lemma in three steps. 
	
	{\it Step 1}. Observability inequality for the dual system.
	
	Let $\psi^2=\sum_{k=1}^n \psi_k^2(x_2)X_k(x_1)\in \mcE_n$. For each $k\in\{1,\dots,n\}$, let $\psi_k$ be the solution to the backward problem
	\begin{equation*}
		\begin{cases} 
			-\pt_t \psi_k - \pt_{x_2} \left(x_2^\al \pt_{x_2} \psi_k\right) + \mu_k x_2^\al \psi_k = 0 & \text{in }\La\ts (t_1,t_2), \\ 
			\psi_k(0,t) = \psi_k(1,t) = 0 & \text{for } t\in (t_1,t_2), \\ 
			\psi_k(x_2,t_2) = \psi_k^2(x_2) & \text{in }\La,
		\end{cases}
	\end{equation*} 
	and set $\psi(x_1,x_2,t):=\sum_{k=1}^n \psi_k(x_2,t)X_k(x_1)$. Applying Corollary \ref{11.27.C1} to $\vp_k(x_2,s):=\psi_k(x_2,t_2-s)$ on $(0,\tau)$ (with $\tau=t_2-t_1\leq T$) and summing over $k=1,\dots,n$, we obtain 
	\begin{equation*}
		\sum_{k=1}^n \|\psi_k(\cdot, t_1)\|_{L^2(\La)}^2\leq C_1e^{C_2\tau^{-2}+C_3n^\be}\int_{t_1}^{t_2}\int_{(a_2,b_2)}\sum_{k=1}^n |\psi_k|^2\df x_2\df t. 
	\end{equation*}
	Combining this with \eqref{11.27.8} and $n\leq n^\be$ (since $\be=\f{3-\al}{2}>1$), we deduce that 
	\begin{equation}\label{11.27.13}
		\|\psi(\cdot, t_1)\|_{L^2(\Om)}^2\leq C_1e^{C_2\tau^{-2}+C_3n^\be} \|\psi\|_{L^2(\om\ts(t_1,t_2))}^2. 
	\end{equation}
	
	{\it Step 2}. Construction of the control. 
	
	Set $\mcR \psi^2:=\psi|_{\om\ts (t_1,t_2)}\in L^2(\om\ts (t_1,t_2))$ and define the linear functional
	\begin{equation*}
		\ell: \Ran \mcR\ra \R, \quad \ell(\mcR\psi^2):=-(\vp^1,\psi(\cdot,t_1))_{L^2(\Om)}, 
	\end{equation*}
	where $\Ran\mcR\s L^2(\om\ts (t_1,t_2))$ denotes the range of $\mcR$. By \eqref{11.27.13}, if $\mcR\psi^2 = 0$, then $\psi(\cdot,t_1)=0$, which implies that $\ell$ is well-defined and bounded on $\Ran\mcR$, satisfying
	\begin{equation*}
		|\ell (\mcR\psi^2)|\leq \sqrt{C_1}e^{\f{1}{2}C_2\tau^{-2}+\f{1}{2}C_3n^\be} \|\vp^1\|_{L^2(\Om)}\|\mcR\psi^2\|_{L^2(\om\ts(t_1,t_2))}. 
	\end{equation*}
	Extending $\ell$ continuously to $\ol{\Ran\mcR}$ and setting it to zero on $(\Ran\mcR)^\perp$, we obtain by the Riesz representation theorem a unique function $f_n\in L^2(\om\ts (t_1,t_2))$ such that 
	\begin{equation}\label{11.27.14}
		\int_{t_1}^{t_2}\int_\om f_n\psi\df x\df t=-(\vp^1, \psi(\cdot, t_1))_{L^2(\Om)}, \text{ for all }\psi^2\in \mcE_n, 
	\end{equation}
	and \eqref{11.27.10} holds. Multiplying \eqref{11.15.21} by $\psi_k$ on $(t_1,t_2)$ with $m=k$ and 
	\begin{equation*}
		f_k(\cdot,t)=\int_\La f_n(x_1,\cdot,t)\chi_\om(x_1,\cdot)X_k(x_1)\df x_1, 
	\end{equation*}
	summing over $k=1,\dots,n$, and integrating by parts, we obtain 
	\begin{equation*}
		(\vp(\cdot,t_2),\psi^2)_{L^2(\Om)}=(\vp(\cdot,t_1),\psi(\cdot,t_1))_{L^2(\Om)}+\int_{t_1}^{t_2}\int_\om f_n\psi\df x\df t. 
	\end{equation*}
	Together with \eqref{11.27.14}, this yields 
	\begin{equation*}
		(\vp(\cdot,t_2),\psi^2)_{L^2(\Om)}=0, \text{ for all }\psi^2\in \mcE_n,
	\end{equation*}
	which shows that \eqref{11.27.12} holds.
	
	{\it Step 3}. Energy estimate at time $t_2$.
	
	By the standard energy estimate (or Duhamel's formula) for \eqref{11.15.1} on $(t_1,t_2)$, we have
	\begin{equation*}
		\|\vp(\cdot, t_2)\|_{L^2(\Om)}\leq \|\vp^1\|_{L^2(\Om)}+\sqrt{\tau}\|f_n\|_{L^2(\om\ts(t_1,t_2))}. 
	\end{equation*}
	Squaring both sides, using $\tau\leq T$, and combining the result with \eqref{11.27.10}, we obtain \eqref{11.27.11}. This completes the proof. 
\end{proof}

\begin{theorem}\label{11.27.T2}
	Let $T>0$, $\al\in (0,1)$, and $\om=(a_1,b_1)\ts (a_2,b_2)$ with $0<a_1<b_1<1$ and $0<a_2<b_2<1$. For each $\vp^0\in L^2(\Om)$, there exists a control $f\in L^2(\om\ts (0,T))$ such that the corresponding solution $\vp$ to \eqref{11.15.1} satisfies $\vp(T)=0$. Moreover, 
	\begin{equation}\label{11.27.9}
		\|f\|_{L^2(\om\ts (0,T))}\leq C\|\vp^0\|_{L^2(\Om)}, 
	\end{equation}
	where the positive constant $C$ depends only on $\al, \om$, and $T$. 
\end{theorem}

\begin{proof}
	We prove this theorem in four steps. 
	
	{\it Step 1}. Parameter selection.
	
	We set
	\begin{equation*}
		\rho=\f{1-\al}{4},\quad p=(2-\al)-\rho=\f{7-3\al}{4}. 
	\end{equation*}
	Then
	\begin{equation}\label{11.27.15}
		p-\be=\rho>0,\quad p-2\rho=\f{5-\al}{4}>0. 
	\end{equation}
	Taking 
	\begin{equation*}
		K=T(2^\rho-1),\quad \tau_j=K2^{-j\rho},\quad T_j=T(1-2^{-j\rho}),\quad T_0=0, 
	\end{equation*}
	we have $T_j-T_{j-1}=\tau_j$ for all $j\in\N^*$, and $T_j\uparrow T$ as $j\to\iy$. Setting $T_{j-\f{1}{2}}=T_{j-1}+\f{\tau_j}{2}$, we set
	\begin{equation*}
		M_j=2^{j+j_0} \quad \text{for } j\in\N^*, 
	\end{equation*}
	where $j_0\in\N^*$ is a sufficiently large integer to be specified later. 
	
	{\it Step 2}. Construction of the control on each subinterval.
	
	Set $\vp_{j-1}:=\vp(\cdot, T_{j-1})$ for $j\in\N^*$. On $[T_{j-1},T_{j-\f{1}{2}})$, by Lemma \ref{11.27.L2} (applied with $\tau=\f{\tau_j}{2}$ and $n=M_j$), there exists a control $f_j\in L^2(\om\ts (T_{j-1},T_{j-\f{1}{2}}))$ such that 
	\begin{equation}\label{11.27.18}
		P_{M_j}\vp(\cdot, T_{j-\f{1}{2}})=0, 
	\end{equation}
	and 
	\begin{alignat}{3}
		\|f_j\|_{L^2(\om\ts(T_{j-1},T_{j-\f{1}{2}}))}^2
		&\leq C_1e^{4C_2\tau_j^{-2}+C_3M_j^\be}\|\vp_{j-1}\|_{L^2(\Om)}^2, \label{11.27.16}\\
		\|\vp(\cdot, T_{j-\f{1}{2}})\|_{L^2(\Om)}^2
		&\leq 2(1+C_1T)e^{4C_2\tau_j^{-2}+C_3M_j^\be}\|\vp_{j-1}\|_{L^2(\Om)}^2.  \label{11.27.17}
	\end{alignat}
	On $[T_{j-\f{1}{2}}, T_j)$, we take the control to be zero. From Lemma \ref{11.15.L4} and \eqref{11.27.18}, we obtain 
	\begin{equation}\label{11.27.19}
		\|\vp_j\|_{L^2(\Om)}^2\leq e^{-C_\al \tau_jM_j^{2-\al}}\|\vp(\cdot, T_{j-\f{1}{2}})\|_{L^2(\Om)}^2 \quad \text{with } C_\al=\f{(1-\al)^2\pi^{2-\al}}{4}. 
	\end{equation}
	
	{\it Step 3}. Exponential decay of the state.
	
	Combining \eqref{11.27.17} and \eqref{11.27.19}, for all $j\in\N^*$, we obtain
	\begin{equation*}
		\|\vp_j\|_{L^2(\Om)}^2\leq 2(1+C_1T)e^{4C_2K^{-2}2^{2j\rho }+C_32^{j_0\be}2^{j\be}-C_\al K2^{(2-\al)j_0}2^{jp}}\|\vp_{j-1}\|_{L^2(\Om)}^2. 
	\end{equation*}
	Note that $2-\al>\be$. Choosing $j_0\in\N^*$ sufficiently large such that 
	\begin{equation*}
		C_\al K2^{(2-\al)j_0}\geq 4C_2K^{-2}+C_32^{j_0\be}+\log\left[2(1+C_1T)\right]+1, 
	\end{equation*}
	and using \eqref{11.27.15} (which gives $1 < 2^{2j\rho} < 2^{jp}$ and $2^{j\be} < 2^{jp}$ for all $j\in\N^*$), we deduce that 
	\begin{equation*}
		\|\vp_j\|_{L^2(\Om)}^2\leq e^{-2^{pj}}\|\vp_{j-1}\|_{L^2(\Om)}^2, \text{ for all } j\in\N^*.
	\end{equation*}
	Iterating this inequality for $j\in\N^*$, we get
	\begin{equation}\label{11.27.20}
		\|\vp_j\|_{L^2(\Om)}^2\leq e^{-\sum_{k=1}^j 2^{pk}}\|\vp^0\|_{L^2(\Om)}^2\leq e^{-2^{pj}}\|\vp^0\|_{L^2(\Om)}^2.
	\end{equation}
	
	{\it Step 4}. Construction of the global control.
	
	Define
	\begin{equation*}
		f(x,t)=\sum_{j=1}^\iy \chi_{[T_{j-1},T_{j-\f{1}{2}})}(t)f_j(x,t), \text{ for } (x,t)\in\om\ts(0,T).
	\end{equation*}
	From \eqref{11.27.16} and \eqref{11.27.20}, we deduce that 
	\begin{equation*}
		\|f\|_{L^2(\om\ts (0,T))}^2\leq C_1\left(e^{4C_2K^{-2}2^{2\rho}+C_32^{j_0\be}2^{\be}}+\sum_{j=2}^\iy e^{4C_2K^{-2}2^{2j\rho}+C_32^{j_0\be}2^{j\be}-2^{p(j-1)}}\right)\|\vp^0\|_{L^2(\Om)}^2. 
	\end{equation*}
	Since $p>\max\{2\rho,\be\}$, the negative term $-2^{p(j-1)}$ dominates the positive terms $4C_2K^{-2}2^{2j\rho}+C_32^{j_0\be}2^{j\be}$ as $j\to\iy$, so the series converges. This proves that $f\in L^2(\om\ts (0,T))$ and \eqref{11.27.9} holds. Moreover, since $\vp\in C([0,T];L^2(\Om))$ and $T_j\uparrow T$, passing to the limit as $j\to\iy$ in \eqref{11.27.20} yields $\vp(\cdot,T)=0$. This completes the proof. 
\end{proof}

\begin{remark}
	The construction in Theorem~\ref{11.27.T2} follows the Fourier-mode partial-control and dissipation strategy used in \cite[Section 4.3]{Beauchard}. For the present boundary-degenerate diffusion, Corollary~\ref{11.27.C1} and  Lemma~\ref{11.27.L1} yield a partial control cost of the form $e^{C(1+\tau^{-2}+n^\beta)}$, whereas the uncontrolled high-frequency modes decay at a rate of order $n^{2-\alpha}$. The strict inequality $\beta=\f{3-\al}{2}<2-\alpha$ permits the above iteration.
\end{remark}

\end{document}